\documentclass[10pt]{amsart}

\usepackage{a4wide,graphicx,hyperref}

\newtheorem{theorem}{Theorem}[section]
\newtheorem{MainResult}{Result}
\newtheorem{lemma}[theorem]{Lemma}
\newtheorem{proposition}[theorem]{Proposition}
\newtheorem{corollary}[theorem]{Corollary}
\theoremstyle{definition}
\newtheorem{definition}[theorem]{Definition}
\theoremstyle{remark}
\newtheorem{remark}[theorem]{Remark}

\newtheorem{problem}[theorem]{Problem}
\newtheorem{conjecture}[theorem]{Conjecture}
\newtheorem{example}[theorem]{Example}

\numberwithin{equation}{section}

\title{Fractal tiles associated with shift radix systems}

\author[V.~Berth\'e]{Val\'erie Berth\'e}
\address{LIRMM, CNRS UMR 5506, Universit\'e Montpellier II, 161 rue Ada, 34392 Montpellier Cedex 5, FRANCE}
\email{berthe@lirmm.fr}
\author[A.~Siegel]{Anne Siegel}
\address{IRISA, Campus de Beaulieu, 35042 Rennes Cedex, FRANCE}
\email{Anne.Siegel@irisa.fr}
\author[W.~Steiner]{Wolfgang Steiner}
\address{LIAFA, CNRS UMR 7089, Universit\'e Paris Diderot -- Paris 7,
Case 7014, 75205 Paris Cedex 13, FRANCE}
\email{steiner@liafa.jussieu.fr}
\author[P.~Surer]{Paul Surer}
\address{Chair of Mathematics and Statistics, University of Leoben, A-8700 Leoben, AUSTRIA}
\email{me@palovsky.com}
\author[J. M. Thuswaldner]{J\"org M. Thuswaldner}
\address{Chair of Mathematics and Statistics, University of Leoben, A-8700 Leoben, AUSTRIA}
\email{joerg.thuswaldner@mu-leoben.at}
\thanks{This research was supported by the Austrian Science Foundation (FWF),
project S9610, which is part of the national research network
FWF--S96 ``Analytic combinatorics and probabilistic number theory'',
by the Agence Nationale de la Recherche, grant ANR--06--JCJC-0073
``DyCoNum'', by the ``Amad\'ee'' grant FR--13--2008 and the ``PHC Amadeus'' grant 17111UB}

\date{\today}



\begin{document}

\begin{abstract}
Shift radix systems form a collection of dynamical systems depending on a parameter $\mathbf{r}$ which varies in the $d$-dimensional real vector space. They generalize well-known numeration systems such as beta-expansions, expansions with respect to rational bases, and canonical number systems. Beta-numeration and canonical number systems are known to be intimately related to fractal shapes, such as the classical Rauzy fractal and the twin dragon. These fractals turned out to be important for studying properties of expansions in several settings.

In the present paper we associate a collection of fractal tiles with shift radix systems. We show that for certain classes of parameters $\mathbf{r}$ these tiles coincide with affine copies of the well-known tiles associated with beta-expansions and canonical number systems. On the other hand, these tiles provide natural families of tiles for beta-expansions with (non-unit) Pisot numbers as well as canonical number systems with (non-monic) expanding polynomials.

We also prove basic properties for tiles associated with shift radix systems. Indeed, we prove that under some algebraic conditions on the parameter $\mathbf{r}$ of the shift radix system, these tiles provide multiple tilings and even tilings of the $d$-dimensional real vector space. These tilings turn out to have a more complicated structure than the tilings arising from the known number systems mentioned above. Such a tiling may consist of tiles having infinitely many different shapes. Moreover, the tiles need not be self-affine (or graph directed self-affine).
\end{abstract}

\maketitle

\section{Introduction}

\subsection*{Number systems, dynamics and fractal geometry}
In the last decades, dynamical systems and fractal geometry have
been proved to be deeply related to the study of number systems (see
\emph{e.g.} the survey~\cite{Barat-Berthe-Liardet-Thuswaldner:06}).
Famous examples of fractal tiles that stem from number systems
are given by the twin dragon fractal (upper left part of
Figure~\ref{fracs} below) which is related to expansions of Gaussian
integers in base $-1+i$ (see~\cite[p.~206]{Knuth:98}), or by the
classical Rauzy fractal (upper right part of Figure~\ref{fracs}) which is
related to beta-expansions with respect to the Tribonacci number~$\beta$
(satisfying $\beta^3=\beta^2+\beta+1$; \emph{cf.}~\cite{Rauzy:82,Thurston:89}).
Moreover, we mention  the Hokkaido tile
(lower left part of Figure~\ref{fracs}) which is related to the
smallest Pisot number (see~\cite{Akiyama:02}) and has been studied
frequently in the literature.

For several notions of number system geometric and dynamical
considerations on fractals imply various non-trivial number
theoretical properties. The boundary of these fractals is intimately
related to the addition of~1 in the underlying number system.
Moreover, the fact
that the origin is an inner point of such a fractal has several
implications. For instance, in beta-numeration as well as for the
case of canonical number systems it implies that the underlying
number system admits finite expansions (all these relations are
discussed in \cite{Barat-Berthe-Liardet-Thuswaldner:06}). In the
case of the classical Rauzy fractal, this allows the computation of
best rational simultaneous approximations of the vector
$(1,1/\beta,1/\beta^2)$, where $\beta$ is the Tribonacci number (see
\cite{CHM01,HM06b}). Another example providing a relation between
fractals and numeration is given by the local spiral shape of the
boundary of the Hokkaido tile: by constructing a realization of the
natural extension for the so-called beta-transformation, one proves
that unexpected non-uniformity phenomena appear in the
beta-numeration associated with the smallest Pisot number. Indeed,
the smallest positive real number that can be approximated by rational
numbers with non-purely periodic beta-expansion is an irrational
number slightly smaller than $2/3$ (see \cite{AFSS,AS:05}).

All the fractals mentioned so far are examples of the new class of
tiles associated with shift radix systems which forms the main object
of the present paper.

\begin{figure}[h]
\centering \leavevmode
\includegraphics[width=0.28\textwidth]{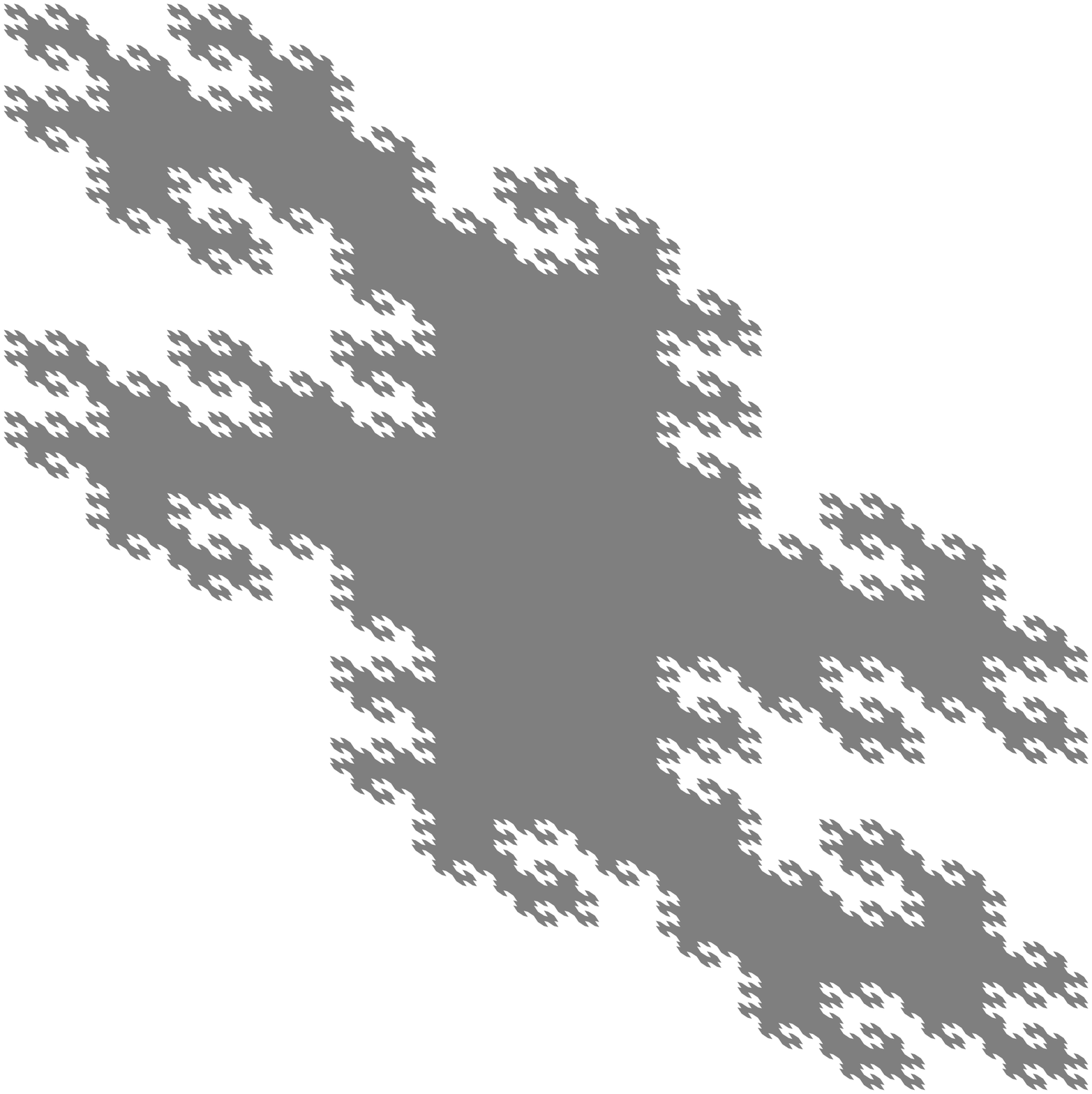} \hskip 2cm
\includegraphics[width=0.28\textwidth]{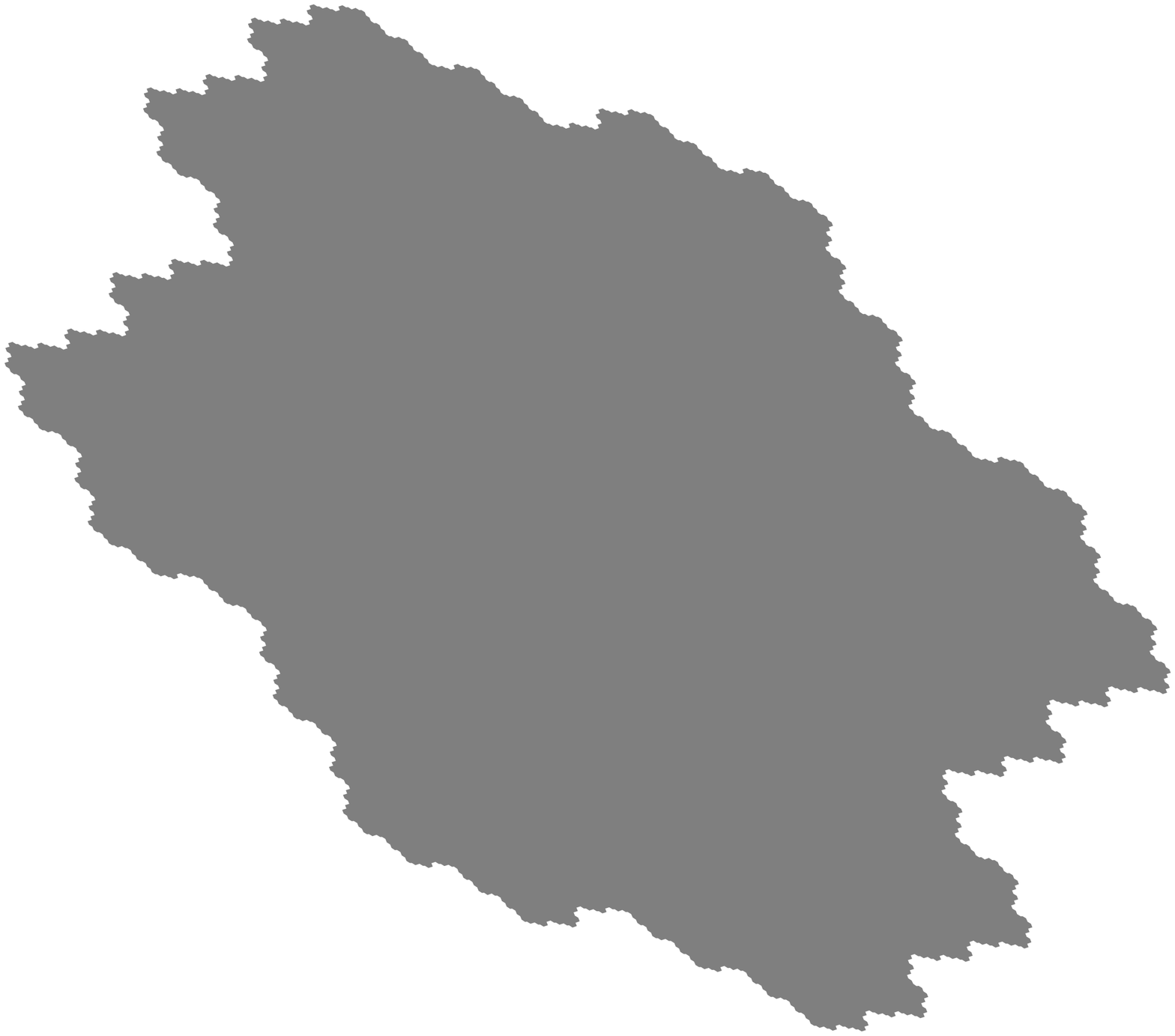}\\
\includegraphics[width=0.28\textwidth]{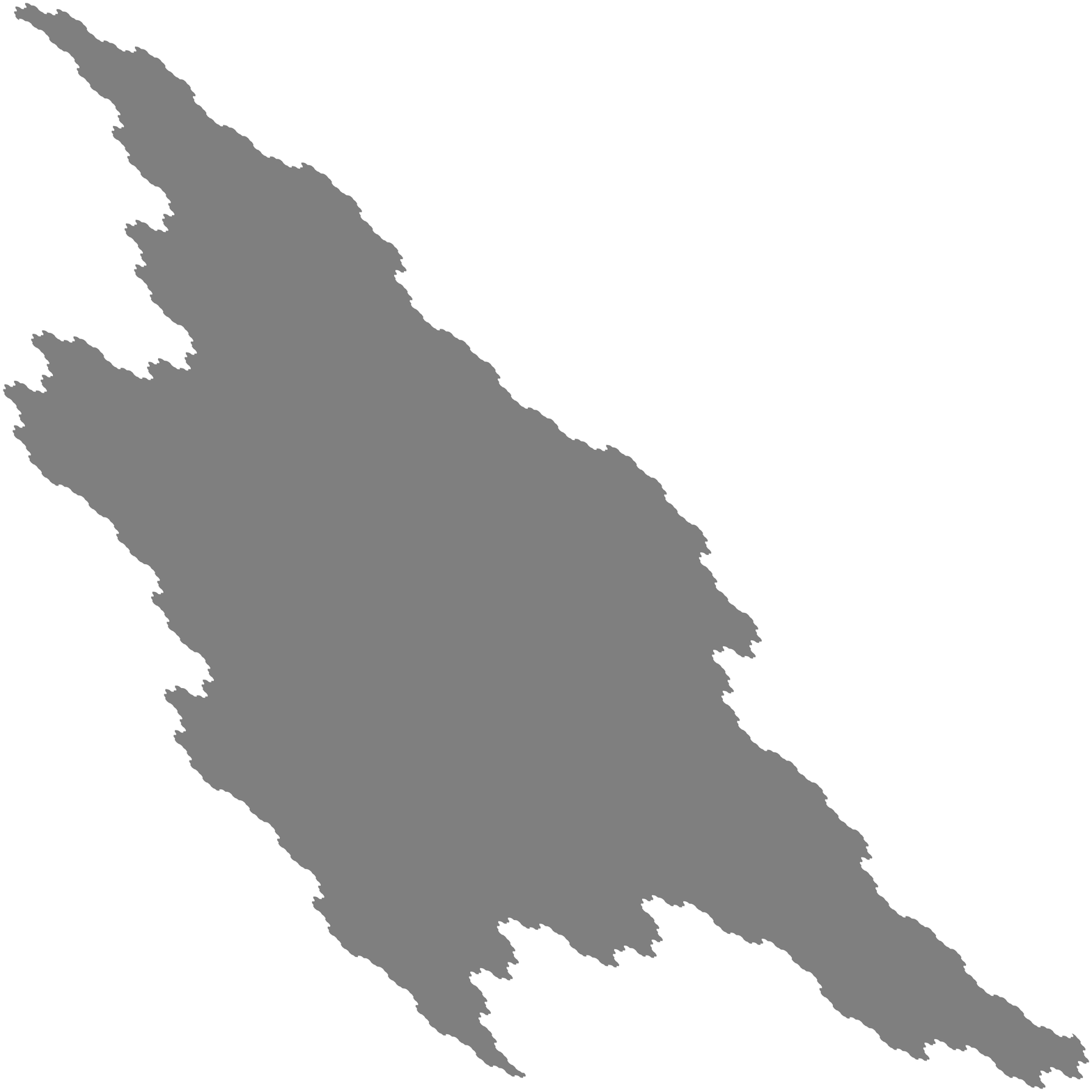}\hskip 2cm
\includegraphics[width=0.28\textwidth]{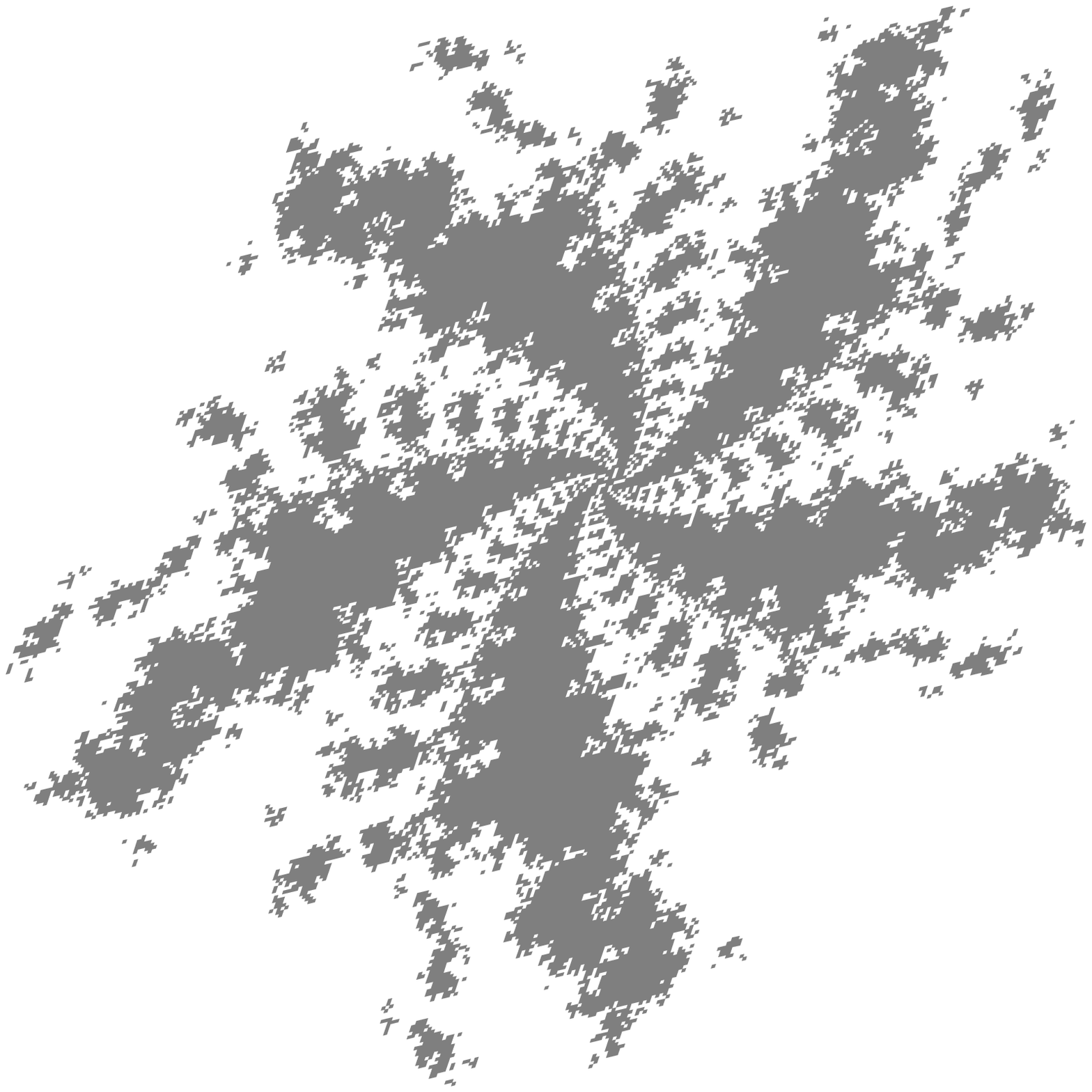}
\caption{Examples of SRS tiles: The upper left tile is the so-called
twin dragon fractal (\emph{cf.} \cite{Knuth:98}), right beside it is
the well-known Rauzy fractal associated with  the Tribonacci number
(\emph{cf.} \cite{Rauzy:82}). The lower left tile is known as
``Hokkaido fractal'' and corresponds to the smallest Pisot number
which has minimal polynomial $x^3-x-1$ (\emph{cf.}
\cite{Akiyama:02}). The lower right one seems to be new, and is an
SRS tile associated with the parameter $\mathbf{r}=(9/10,-11/20)$.}
\label{fracs}
\end{figure}

\subsection*{Shift radix systems: a common dynamical formalism}
Shift radix systems have been proposed in
\cite{Akiyama-Borbeli-Brunotte-Pethoe-Thuswaldner:05} to unify
various notions of radix expansion such as beta-expansions (see
\cite{Frougny-Solomyak:92,Parry:60,Renyi:57}), canonical number
systems (CNS for short, see \cite{Kovacs-Pethoe:91,Petho:91}) and
number systems with respect to rational bases (in the sense of
\cite{Akiyama-Frougny-Sakarovitch:07}) under a common roof. Instead
of starting with a base number (or base polynomial), one considers a
vector $\mathbf{r} \in \mathbb{R}^d$ and defines the mapping
$\tau_\mathbf{r}:\mathbb{Z}^d \to \mathbb{Z}^d$ by
\[
\tau_\mathbf{r}(\mathbf{z})=
(z_1,\ldots,z_{d-1},-\lfloor\mathbf{r}\mathbf{z}\rfloor)^t \quad \big(\mathbf{z}=(z_0,\ldots,z_{d-1})^t\big).
\]
Here, $\mathbf{r}\mathbf{z}$ denotes the scalar product of the vectors $\mathbf{r}$ and $\mathbf{z}$; moreover, for $x$ being a real number, $\lfloor x \rfloor$ denotes the largest integer less than or equal to $x$.
We call $(\mathbb{Z}^d, \tau_\mathbf{r})$ a \emph{shift radix system} (SRS, for short).
A~vector $\mathbf{r}$ gives rise to an SRS with \emph{finiteness property} if each
integer vector $\mathbf{z}\in\mathbb{Z}^d$ can be finitely expanded
with respect to the vector $\mathbf{r}$, that is, if for each
$\mathbf{z}\in\mathbb{Z}^d$ there is an $n\in\mathbb{N}$ such that
the $n$-th iterate of $\tau_\mathbf{r}$ satisfies
$\tau_\mathbf{r}^n(\mathbf{z})=\mathbf{0}$.

In the papers written on SRS so far (see \emph{e.g.}\
\cite{Akiyama-Borbeli-Brunotte-Pethoe-Thuswaldner:05,Akiyama-Brunotte-Pethoe-Thuswaldner:06,Akiyama-Brunotte-Pethoe-Thuswaldner:07,Akiyama-Brunotte-Pethoe-Thuswaldner:08,Surer:07}),
relations between SRS and well-known notions of number system such
as beta-expansions with respect to unit Pisot numbers and canonical
number systems with a monic expanding polynomial\footnote{An \emph{expanding polynomial} is a polynomial each of whose roots has modulus strictly larger than one.} have been
established (we will give a short account on these relations in
Section~\ref{cns} and~\ref{beta}, respectively). In particular, SRS
turned out to be a fruitful tool in order to deal with the problem
of finite representations in these number systems.

\subsection*{SRS tiles: an extension of  arithmetically meaningful known constructions}
In the present paper, we study geometric properties of SRS: we introduce a family of tiles for each $\mathbf{r} \in \mathbb{R}^d$ with contractive and regular companion matrix, and
prove geometric properties of these tiles. Figure~\ref{fracs} shows four examples of such tiles, called \emph{SRS tiles}.

We prove that our definition unifies the notions of self-affine
tiles known for CNS with respect to monic polynomial
(see~\cite{Barat-Berthe-Liardet-Thuswaldner:06,Katai-Koernyei:92})
and beta-numeration related to a unit Pisot number
(see~\cite{Akiyama:02}).

\begin{MainResult}
The following relations between SRS tiles and classical tiles related to number systems hold:
\begin{itemize}
\item If the parameter ${\bf r}$ of an SRS is related to a monic expanding polynomial over ${\mathbb Z}$, the SRS tile is a linear image of the self-affine
tile associated with this polynomial.
\item If the parameter ${\bf r}$ of an SRS comes from a unit Pisot number, the SRS tile is a linear image of the central tile associated with the corresponding beta-numeration.
\end{itemize}
\end{MainResult}

It is well-known that these classical tiles are self-affine and
associated with (multiple) tilings that are highly structured (see
\emph{e.g.}\ the
surveys~\cite{Akiyama-Thuswaldner:04,Berthe-Siegel:05}). Indeed,
these tiles  satisfy a set equation expressed as a graph-directed
iterated function system (GIFS) in the sense of Mauldin and
Williams~\cite{Mauldin-Williams:88} which means that each tile can
be decomposed with respect to the action of an affine mapping into
several copies of a finite number of tiles, with the decomposition
being produced by a finite graph
(see~\cite{Barat-Berthe-Liardet-Thuswaldner:06}).

In the present SRS situation, we prove that our construction gives
rise to new classes of tiles, in particular we want to emphasize on
tiles related to a non-monic expanding polynomial or to  a non-unit
Pisot number. In both cases, we are actually able to extend the
usual definition of tiles used for  unit Pisot numbers and monic
polynomials.  The  tiles defined in the usual way are not
satisfactory in this more general setting, since they always produce
overlaps in their self-affine decomposition. A first strategy to
remove overlaps consists in enlarging the space of representation by
adding arithmetic components ($p$-adic factors) as proposed
in~\cite{ABBS:08}. However, such tiles are of limited topological
importance since they have totally disconnected factors.

We prove that our construction, however, allows to insert arithmetic
criteria in the construction of the tiles: roughly speaking, the SRS
mapping $\tau_\mathbf{r}$ naturally selects points in appropriate
submodules. This arithmetic selection process removes the overlaps.

\begin{MainResult}
SRS provide a  natural collection of  tiles for number systems related to non-monic expanding polynomials  as well as to non-unit Pisot numbers.
\end{MainResult}

Nonetheless, the geometrical structure of these tiles is much more
complicated than the structure of the classical ones. As we shall
illustrate for $\mathbf{r}=(-2/3)$, there may be infinitely many
shapes of tiles associated with certain parameters~$\mathbf{r}$.
Actually, the description of SRS tiles requires a set equation that
cannot be captured by a finite graph. Equation
(\ref{eq:decomposition}) suggests that an infinite hierarchy of set
equations is needed to describe an SRS tile. Therefore, SRS tiles in
general cannot be regarded as GIFS attractors. Furthermore, an SRS
tile is not always equal to the closure of its interior
(Example~\ref{ex:pp} exhibits a case of an SRS tile that is equal to
a single point). Also, due to the lack of a GIFS structure, we have
no information on the measure of the boundary of tiles.

\subsection*{Tiling properties}
Despite of their complicated geometrical structure, we are able to
prove tiling properties for SRS tiles. We first prove that, for each
fixed parameter~$\mathbf{r}$,  the associated tiles form a covering,
extending the results known for unit Pisot numbers and monic
expanding polynomials. In the classical cases, exhibiting a multiple
tiling (that is, a covering with an almost everywhere constant
degree) from this covering is usually done by exploiting specific
features of the tiling together with the GIFS structure of the
tiles: this allows to transfer  local information to the whole
 space in order to obtain  global information.  The finiteness property is then used to prove
that  $\mathbf{0}$ belongs to only one tile, leading to a tiling. In
the present SRS situation, however, this strategy does not work any
more. The finiteness property is still equivalent to the fact that
$\mathbf{0}$ belongs to exactly one  SRS tile, but the fact that SRS
tiles are no longer GIFS attractors prevents us from spreading this
local information to the whole space. In order to get global
information, we impose additional algebraic conditions on the
parameter~$\mathbf{r}$, and we use these conditions to exhibit a dense
set of points that are proved to belong to a fixed number of tiles.
This leads to the multiple tiling property.

Our main result can thus be stated as follows (for the corresponding definitions and  for a more precise statement, see
Section~\ref{sec:multiple}, in particular Theorem~\ref{theo:mdense} and Corollary~\ref{cor:weaktiling}):

\begin{MainResult} Let $\mathbf{r}$ be an SRS parameter with contractive and regular companion matrix, and assume that either $\mathbf{r}\in\mathbb{Q}^d$ or $\mathbf{r}$ is related to a Pisot number or  $\mathbf{r}$ has algebraically independent coordinates. Then the following assertions hold.
\begin{itemize}
\item
The collection of SRS tiles associated with the parameter $\mathbf{r}$ forms a weak multiple tiling.
\item
If the finiteness property is satisfied, then the collection of SRS tiles forms a weak tiling.
\end{itemize}
\end{MainResult}

By weak  tiling, we mean here that the tiles cover the whole space
and have disjoint interiors. Stating a ``strong'' tiling property
would require to have information on the boundaries of the tiles,
which is deeply intricate if $d\ge 2$ since the tiles are no longer
GIFS attractors, and deserves a specific study. For $d=1$ the
situation becomes easier. Indeed, we will prove in Theorem~\ref{d1}
that in this case SRS tiles are (possibly degenerate) intervals
which form a tiling of $\mathbb{R}$.

\subsection*{Structure of the paper}
In Section~\ref{srsrep}, we introduce a way of representing integer
vectors by using the shift radix transformation~$\tau_\mathbf{r}$.
In Section~\ref{sec:SRStiles}, SRS tiles are defined and fundamental
geometric properties of them are studied. Section~\ref{sec:multiple}
is devoted to tiling properties of SRS tiles. We show that SRS tiles
form tilings or multiple tilings for large classes of parameters
$\mathbf{r}\in\mathbb{R}^d$. As the tilings are no longer
self-affine we have to use new methods in our proofs. In
Section~\ref{cns}, we analyse the relation between tiles associated
with expanding polynomials and SRS tiles more closely. We prove that
the SRS tiles given by a monic CNS parameter $\mathbf{r}$ coincide
up to a linear transformation with the self-affine CNS tile. For
non-monic expanding polynomials, SRS tiles give rise to a new class
of tiles. At the end of this section, we prove for the parameter
$\mathbf{r}=(-2/3)$ that the associated tiling has infinitely many
shapes of tiles. In Section~\ref{beta}, after proving that the shift
radix transformation is conjugate to the beta-transformation
restricted to $\mathbb{Z}[\beta]$, we investigate the relation
between beta-tiles and SRS tiles. It turns out that beta-tiles
associated with a unit Pisot number are linear images of SRS tiles
related to a parameter associated with this Pisot number. Moreover,
we define a new class of tiles associated with non-unit Pisot
numbers. The paper ends with a short section in which conjectures and possible directions of future research related to the topic of the present paper are discussed.

\section{The SRS representation}\label{srsrep}

\subsection{Definition of shift radix systems and their parameter domains}
Shift radix systems are  dynamical systems  defined  on
$\mathbb{Z}^d$ as follows
(see~\cite{Akiyama-Borbeli-Brunotte-Pethoe-Thuswaldner:05}).

\begin{definition}[Shift radix system, finiteness property]\label{SRSdef}
For $\mathbf{r}=(r_0,\ldots, r_{d-1}) \in \mathbb{R}^d$,
$d\ge1$, set
\begin{eqnarray*}
\tau_\mathbf{r}:\; \mathbb{Z}^d &\rightarrow& \mathbb{Z}^d,\\
\mathbf{z}=(z_0,z_1,\ldots,z_{d-1})^t &\mapsto& (z_1,\ldots,z_{d-1},-\lfloor\mathbf{r}\mathbf{z}\rfloor)^t,
\end{eqnarray*}
where $\mathbf{r}\mathbf{z}$ denotes the scalar product of
$\mathbf{r}$ and $\mathbf{z}$. We call the dynamical system
$(\mathbb{Z}^d,\tau_\mathbf{r})$ a \emph{shift radix system} (SRS,
for short). The SRS parameter $\mathbf{r}$ is said to be
\emph{reduced} if $r_0 \neq 0$.

We say that $(\mathbb{Z}^d,\tau_\mathbf{r})$ satisfies the
\emph{finiteness property} if for every $\mathbf{z} \in
\mathbb{Z}^d$ there exists some $n \in \mathbb{N}$ such that
$\tau_\mathbf{r}^n(\mathbf{z})=\mathbf{0}$.
\end{definition}

\begin{remark}
If $r_0=0$, then every vector $\tau_\mathbf{r}^n(\mathbf{z})$, $\mathbf{z} \in \mathbb{Z}^d$, $n \ge 1$, can be easily obtained from the SRS  $(\mathbb{Z}^{d-1}, \tau_{\mathbf{r}'})$ with $\mathbf{r}' = (r_1,\ldots,r_{d-1})$ since $\tau_{\mathbf{r}'}\circ\pi = \pi\circ\tau_\mathbf{r}$, where $\pi:\, \mathbb{Z}^d \to \mathbb{Z}^{d-1}$ denotes the projection defined by  $\pi(z_0,z_1,\ldots,z_{d-1}) = (z_1,\ldots,z_{d-1})$.
\end{remark}

The companion matrix of $\mathbf{r}=(r_0,\ldots,r_{d-1})$ is denoted by
\[
M_\mathbf{r} := \begin{pmatrix}
0 & 1 & 0 & \cdots & 0 \\
\vdots & 0 & \ddots & \ddots & \vdots \\
\vdots & \vdots & \ddots & \ddots & 0 \\
0 & 0 & \cdots & 0 & 1 \\
-r_{0} & -r_{1} & \cdots & -r_{d-2} &-r_{d-1}
\end{pmatrix} \in \mathbb{R}^{d \times d}.
\]
Its characteristic polynomial is given by $X^d + r_{d-1}X^{d-1}+\cdots+r_1X+r_0$. The matrix $M_\mathbf{r}$ is regular if $r_0\not=0$.  We will work with  $M_\mathbf{r}^{-1}$  (see \emph{e.g.}\  Proposition \ref{ex:tau}), hence   we will assume
${\bf r}$ reduced in all that follows.  For
$\mathbf{z}=(z_0,\ldots,z_{d-1})^t$, note that
\[
M_\mathbf{r}\mathbf{z}=(z_1,\ldots,z_{d-1},-\mathbf{r}\mathbf{z})^t
\]
and thus
\begin{equation} \label{eq:tauR}
\tau_\mathbf{r}(\mathbf{z})=M_\mathbf{r}\mathbf{z}+(0,\ldots,0,\{\mathbf{r}\mathbf{z}\})^t,
\end{equation}
where $\{x\}=x-\lfloor x\rfloor$ denotes the fractional part of $x$.

The sets
\begin{align*}
\mathcal{D}_d & := \big\{\mathbf{r} \in \mathbb{R}^d \mid (\tau_\mathbf{r}^n(\mathbf{z}))_{n \in \mathbb{N}}\ \mbox{is eventually periodic for all}\ \mathbf{z} \in \mathbb{Z}^d\big\}\quad \mbox{and} \\
\mathcal{D}_d^{(0)} & := \big\{\mathbf{r} \in \mathbb{R}^d \mid \tau_\mathbf{r}\ \mbox{is an SRS with finiteness property}\big\}
\end{align*}
are intimately related to the notion of SRS.
Obviously we have $\mathcal{D}_d^{(0)} \subset \mathcal{D}_d$.
Apart from its boundary, the set $\mathcal{D}_d$ is easy to describe.

\begin{lemma}[see \cite{Akiyama-Borbeli-Brunotte-Pethoe-Thuswaldner:05}]
 \label{lem:int}
A point $\mathbf{r} \in \mathbb{R}^d$ is contained in the interior of $\mathcal{D}_d$ if and only if the spectral radius of the companion matrix  $M_\mathbf{r}$ generated by $\mathbf{r}$ is strictly less than~$1$.
\end{lemma}

\begin{figure}[h]
\centering \leavevmode
\includegraphics[width=0.25\textwidth]{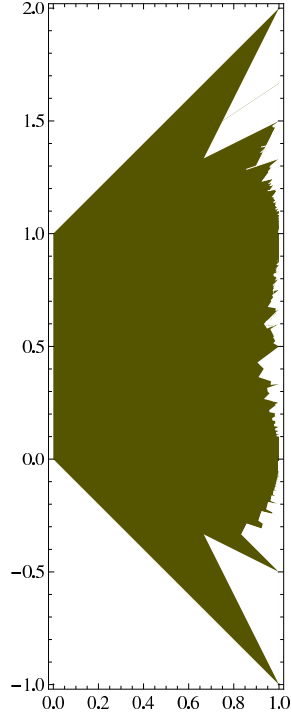} \hskip 2cm
\caption{An approximation of $\mathcal{D}_2^{(0)}$.} \label{srs2}
\end{figure}

The shape of $\mathcal{D}_d^{(0)}$ is of a more complicated nature.
While $\mathcal{D}_1^{(0)}=[0,1)$ is quite easy to describe
(\emph{cf.}\
\cite[Proposition~4.4]{Akiyama-Borbeli-Brunotte-Pethoe-Thuswaldner:05}),
already for $d=2$ no complete description of $\mathcal{D}_2^{(0)}$
is known. The available partial results
(see~\cite{Akiyama-Brunotte-Pethoe-Thuswaldner:06,Surer:07})
indicate that $\mathcal{D}_d^{(0)}$ is of a quite irregular
structure for $d\ge 2$. Nonetheless, an algorithm is given in
\cite[Proposition~4.4]{Akiyama-Borbeli-Brunotte-Pethoe-Thuswaldner:05}
which decides whether a given $\mathbf{r}$ belongs to
$\mathcal{D}_d^{(0)}$. This algorithm was used to draw the
approximation of $\mathcal{D}_2^{(0)}$ depicted in
Figure~\ref{srs2}.

Apart from the easier case $d=1$, in this paper, we are going to
consider three classes of   points  ${\mathbf r}$ assumed  to be reduced which belong to $\mathcal{D}_d$.
The first two classes are dense in $\mathrm{int}(\mathcal{D}_d)$,
while the third one has full measure in
$\mathrm{int}(\mathcal{D}_d)$.
\begin{enumerate}
\item
$\mathbf{r} \in \mathbb{Q}^d \cap \mathrm{int}(\mathcal{D}_d)$.
This class
includes the parameters
$\mathbf{r}=\big(\frac{1}{a_0},\frac{a_{d-1}}{a_0},\ldots,\frac{a_1}{a_0}\big)
\in \mathrm{int}(\mathcal{D}_d)$ with $a_0,\ldots,a_{d-1}\in \mathbb{Z}$, which correspond to expansions with respect to monic expanding polynomials, including CNS. If the first
coordinate of $\mathbf{r}$ has numerator greater than one, this extends to non-monic polynomials in a natural way (see Section~\ref{cns} for details).
\item
$\mathbf{r}=(r_0,\ldots,r_{d-1})$ is obtained by decomposing the
minimal polynomial of a Pisot number~$\beta$ as
$(x-\beta)(x^d+r_{d-1}x^{d-1}+\cdots+r_0)$. In view of
Lemma~\ref{lem:int}, this implies that $\mathbf{r}\in
\mathrm{int}(\mathcal{D}_d)$, and by
\cite{Akiyama-Brunotte-Pethoe-Thuswaldner:08} this set of parameters
is dense in $\mathrm{int}(\mathcal{D}_d)$. These parameters
correspond to beta-numeration with respect to Pisot numbers (see
Section~\ref{beta}). Note that even non-unit Pisot numbers are
covered here.
\item
$\mathbf{r}=(r_0,\ldots,r_{d-1})\in \mathrm{int}(\mathcal{D}_d)$ with algebraically independent coordinates $r_0,\ldots,r_{d-1}$.
\end{enumerate}

\subsection{SRS representation  of  $d$-dimensional integer vectors}  \label{subsec:exp}
For $\mathbf{r} \in \mathbb{R}^d$ we can use the SRS transformation
$\tau_\mathbf{r}$ to define an expansion for $d$-dimensional integer
vectors.

\begin{definition}[SRS representation]\label{srsex}
Let $\mathbf{r} \in \mathbb{R}^d$.
For $\mathbf{z} \in \mathbb{Z}^d$, the \emph{SRS representation} of $\mathbf{z}$ with respect to $\mathbf{r}$ is defined to be the sequence $(v_1,v_2,v_3,\ldots)$, with $v_n=\big\{\mathbf{r}\tau_\mathbf{r}^{n-1}(\mathbf{z})\big\}$ for all $n\ge1$.

The representation is said to be \emph{finite} if there is an $n_0$
such that $v_n=0$ for all $n\ge n_0$. It is said to be
\emph{eventually periodic} if there are $n_0$, $p$ such that
$v_n=v_{n+p}$ for all $n\ge n_0$.
\end{definition}

By definition, every SRS representation is finite if
$\mathbf{r}\in\mathcal{D}_d^{(0)}$, and every SRS representation is
eventually periodic if $\mathbf{r}\in\mathcal{D}_d$. The following
simple properties of SRS representations show that integer vectors
can be expanded according to $\tau_\mathbf{r}$. The matrix
$M_\mathbf{r}$ acts as a base and the vectors $(0,\ldots,0,v_j)$ are
the digits.

\begin{lemma}\label{lum}
Let $\mathbf{r} \in \mathbb{R}^d$ and $(v_1,v_2,\ldots)$ be the SRS representation of $\mathbf{z} \in \mathbb{Z}^d$ with respect to~$\mathbf{r}$.
Then the following properties hold for all $n \in \mathbb{N}$:
\begin{enumerate}
\item \label{lum1}
$0\le v_n<1$,
\item \label{lum2}
$\tau_\mathbf{r}^n(\mathbf{z})$ has the SRS representation $(v_{n+1},v_{n+2},v_{n+3},\ldots)$,
\item \label{lum3}
we have
\begin{equation} \label{eq:tauk}
M_\mathbf{r}^n \mathbf{z} = \tau_\mathbf{r}^n(\mathbf{z}) - \sum_{j=1}^n M_\mathbf{r}^{n-j} (0,\ldots,0,v_j)^t.
\end{equation}
\end{enumerate}
\end{lemma}

\begin{proof}
Assertions (1) and (2) follow immediately from the definition of the SRS representation.
By iterating (\ref{eq:tauR}), we obtain
\[
\tau_\mathbf{r}^n(\mathbf{z}) = M_\mathbf{r}^n\mathbf{z} + \sum_{j=1}^n M_\mathbf{r}^{n-j}\big(0,\ldots,0,\big\{\mathbf{r}\tau_\mathbf{r}^{j-1}(\mathbf{z})\big\}\big)^t,
\]
which yields (\ref{eq:tauk}).
\end{proof}

Note that the set of possible SRS digits $(0,\ldots,0,v)$ is
infinite unless $\mathbf{r} \in \mathbb{Q}^d$.

We prove that the SRS representation is unique in the following sense.

\begin{proposition}\label{ex:tau}
Let $\mathbf{r} \in \mathbb{R}^d$ be reduced and suppose that the
SRS-representation of an element $\mathbf{z}_0 \in \mathbb{Z}^d$ is
$(v_1,v_2,v_3,\ldots)$. Assume that for some reals
$v_0,v_{-1},\ldots,v_{-n+1} \in [0,1)$, $n\in\mathbb{N}$, we have
\[
\mathbf{z}_{-k}:=M_\mathbf{r}^{-k}\Big(\mathbf{z}_0-\sum_{j=0}^{k-1}
M_\mathbf{r}^j(0,\ldots,0,v_{-j})^t\Big) \in \mathbb{Z}^d
\quad\mbox{for all}\ 1 \leq k \leq n.
\]
Then $\tau_\mathbf{r}^n(\mathbf{z}_{-n})=\mathbf{z}_0$ and
$\mathbf{z}_{-n}$ has the SRS representation
$(v_{-n+1},v_{-n+2},v_{-n+3},\ldots)$.
\end{proposition}

\begin{proof}
The assertion is obviously true for $n=0$.
Now continue by induction on $n$.
We have
\[
\tau_\mathbf{r}(\mathbf{z}_{-n}) =
\tau_\mathbf{r}\big(M_\mathbf{r}^{-1}\big(\mathbf{z}_{-n+1}-(0,\ldots,0,v_{-n+1})^t\big)\big)
= \mathbf{z}_{-n+1} - (0,\ldots,0,v_{-n+1})^t +
(0,\ldots,0,\{\mathbf{r}\mathbf{z}_{-n}\})^t.
\]
Since $\tau_\mathbf{r}(\mathbf{z}_{-n})\in\mathbb{Z}^d$, we obtain
$\tau_\mathbf{r}(\mathbf{z}_{-n})=\mathbf{z}_{-n+1}$ and
$v_{-n+1}=\{\mathbf{r}\mathbf{z}_{-n}\}$. Therefore the first SRS
digit of $\mathbf{z}_{-n}$ is $v_{-n+1}$. By induction, we have
$\tau_\mathbf{r}^n(\mathbf{z}_{-n})=\mathbf{z}_0$, hence, the SRS
representation of $\mathbf{z}_{-n}$ is
$(v_{-n+1},v_{-n+2},v_{-n+3},\ldots)$.
\end{proof}

\section{Definition and first properties of SRS tiles} \label{sec:SRStiles}

We now define a new type of tiles based on the mapping
$\tau_\mathbf{r}$ when the matrix $M_\mathbf{r}$ is contractive and $\mathbf{r}$ is reduced. By
analogy with the definition of tiles for other dynamical systems
(see \emph{e.g.}
\cite{Akiyama:98,Berthe-Siegel:05,Katai-Koernyei:92,Rauzy:82,Scheicher-Thuswaldner:01,Thurston:89}),
we consider elements of $\mathbb{Z}^{d}$ which are mapped to a given
$\mathbf x\in \mathbb{Z}^d$ by $\tau_\mathbf{r}^n$, renormalize by a
multiplication with $M_\mathbf{r}^n$, and let $n$ tend to $\infty$.
To build this set, we thus consider vectors whose SRS expansion
coincides with the expansion of $\mathbf{x}$ up to an added finite
prefix and we then renormalize this expansion. We will see in
Sections~\ref{cns} and~\ref{beta} that some of these tiles
are related to well-known types of tiles, namely CNS tiles and beta-tiles.
We recall that four examples of central SRS tiles are depicted in
Figure~\ref{fracs}.

\subsection{Definition of SRS tiles} \label{defsrs}

An SRS tile will turn out to be the limit of the sequence of compact
sets $(M_\mathbf{r}^n \tau_\mathbf{r}^{-n}(\mathbf{x}))_{n\ge0}$
with respect to the Hausdorff metric. As it is \emph{a priori} not
clear that this limit exists, we first define the tiles as lower
Hausdorff limits of these sets and then show that the Hausdorff
limit exists. Recall that the lower Hausdorff limit $\mathop{\rm
Li}_{n\to\infty}A_n$ of a sequence $(A_n)$ of subsets of
$\mathbb{R}^d$ is the (closed) set of all $\mathbf{t}\in \mathbb{R}^d$ having
the property that each neighborhood of~$\mathbf{t}$ intersects $A_n$
provided that $n$ is sufficiently large. If the sets $A_n$ are
compact and $(A_n)$ is a Cauchy sequence w.r.t.\ the Hausdorff metric~$\delta$, then the Hausdorff limit $\mathop{\rm
Lim}\limits_{n\to\infty}A_n$ exists and
 \[
\mathop{\rm Lim}\limits_{n\to\infty}A_n = \mathop{\rm
Li}\limits_{n\to\infty}A_n
\]
(see \cite[Chapter~II, \S29]{Kuratowski:66} for details on Hausdorff
limits).

\begin{definition}[SRS tile]\label{srstdef}
Let $\mathbf{r}\in \mathrm{int}(\mathcal{D}_d)$ be reduced and $\mathbf{x} \in \mathbb{Z}^{d}$.
The \emph{SRS tile} associated with~$\mathbf{r}$ is defined as
\[
\mathcal{T}_\mathbf{r}(\mathbf{x}) := \mathop{\rm
Li}_{n\rightarrow\infty} M_\mathbf{r}^n
\tau_\mathbf{r}^{-n}(\mathbf{x}).
\]
\end{definition}

\begin{remark} \label{r:anne}
Note that this means that $\mathbf{t}\in\mathbb{R}^d$ is an element
of $\mathcal{T}_\mathbf{r}(\mathbf{x})$ if and only if there exist
vectors $\mathbf{z}_{-n}\in\mathbb{Z}^d$, $n\in\mathbb{N}$, such
that $\tau_\mathbf{r}^n(\mathbf{z}_{-n})=\mathbf{x}$ for all
$n\in\mathbb{N}$ and
$\lim_{n\to\infty}M_\mathbf{r}^n\mathbf{z}_{-n}=\mathbf{t}$.
\end{remark}

We will see in Theorem~\ref{seteq} that the lower limit in
Definition~\ref{srstdef} is equal to the limit with respect to the
Hausdorff metric.

\subsection{Compactness, a set equation, and a covering property of $\mathcal{T}_\mathbf{r}(\mathbf{x})$}

In this subsection we will show that each SRS tile
$\mathcal{T}_\mathbf{r}(\mathbf{x})$ is a compact set that can be
decomposed into subtiles which are obtained by multiplying other
tiles by $M_\mathbf{r}$. Moreover, we prove that for each reduced
$\mathbf{r} \in \mathrm{int}(\mathcal{D}_d)$ the collection
$\{\mathcal{T}_\mathbf{r}(\mathbf{x}) \mid \mathbf{x}\in\mathbb{Z}^d\}$ covers the real vector space $\mathbb{R}^d$.

If $\mathbf{r} \in \mathrm{int}(\mathcal{D}_d)$, then the matrix
$M_\mathbf{r}$ is contractive by Lemma~\ref{lem:int}. Let $\rho<1$
be the spectral radius of $M_\mathbf{r}$, $\rho<\tilde\rho<1$ and
$\|\cdot\|$ be a norm satisfying
\[
\|M_\mathbf{r}\mathbf{x}\|\leq \tilde\rho \|\mathbf{x}\| \quad
\mbox{for all } \mathbf{x} \in \mathbb{R}^d
\]
(for the construction of such a norm see for
instance~\cite[Equation~(3.2)]{Lagarias-Wang:96a}). Then we have
\begin{equation} \label{eq:R}
R := \sum_{n=0}^\infty \big\|M_\mathbf{r}^n(0,\ldots,0,1)^t\big\| \le \frac{\|(0,\ldots,0,1)^t\|}{1-\tilde\rho}\,.
\end{equation}

\begin{lemma}\label{endlich}
Let $\mathbf{r} \in \mathrm{int}(\mathcal{D}_d)$ be reduced.
Then every  $\mathcal{T}_\mathbf{r}(\mathbf{x})$, $\mathbf{x}\in\mathbb{Z}^d$, is contained in the closed ball of radius $R$ with center~$\mathbf{x}$.
Therefore the cardinality of the sets $\{\mathbf{x} \in \mathbb{Z}^d \mid \mathbf{t} \in \mathcal{T}_\mathbf{r}(\mathbf{x})\}$ is uniformly bounded in $\mathbf{t}\in \mathbb{R}^d$.
Furthermore, the family of SRS tiles $\{\mathcal{T}_\mathbf{r}(\mathbf{x}) \mid \mathbf{x}\in\mathbb{Z}^d\}$ is \emph{locally finite}, that is, any open ball meets only a finite number of tiles of the family.
\end{lemma}

\begin{proof}
Let $\mathbf{x}\in\mathbb{Z}^d$,
$\mathbf{t}\in\mathcal{T}_\mathbf{r}(\mathbf{x})$ and
$\mathbf{z}_{-n}$ as in Remark~\ref{r:anne}.
Let the SRS representation of $\mathbf{z}_{-n}$ be
$(v_{-n+1}^{(n)},v_{-n+2}^{(n)},v_{-n+3}^{(n)},\ldots)$. Then by
(\ref{eq:tauk}) we have
\[
M_\mathbf{r}^n\mathbf{z}_{-n} = \mathbf{x} -
\sum_{j=1}^n M_\mathbf{r}^{n-j}(0,\ldots,0,v_{-n+j}^{(n)})^t,
\]
thus $\|M_\mathbf{r}^n\mathbf{z}_{-n}-\mathbf{x}\|<R$ and, hence,
$\|\mathbf{t}-\mathbf{x}\|\le R$. The uniform boundedness of the
cardinality of $\{\mathbf{x} \in \mathbb{Z}^d \mid \mathbf{t} \in
\mathcal{T}_\mathbf{r}(\mathbf{x})\}$ and the local finiteness
follow immediately.
\end{proof}

\begin{lemma} \label{lem:Hausdorff}
Let $\mathbf{r} \in \mathrm{int}(\mathcal{D}_d)$ be reduced and denote by
$\delta(\cdot,\cdot)$ the Hausdorff metric induced by a norm
satisfying (\ref{eq:R}). Then
$(M_\mathbf{r}^n\tau_\mathbf{r}^{-n}(\mathbf{x}))_{n\ge 0}$ is a
Cauchy sequence w.r.t. $\delta$, in particular,
\[
\delta\big(M_\mathbf{r}^n\tau_\mathbf{r}^{-n}(\mathbf{x}),M_\mathbf{r}^{n+1}\tau_\mathbf{r}^{-n-1}(\mathbf{x})\big)
\le \tilde\rho^n\|(0,\ldots,0,1)^t\|.
\]
\end{lemma}

\begin{proof}
Let $\mathbf{t}\in M_\mathbf{r}^n\tau_\mathbf{r}^{-n}(\mathbf{x})$.
Then
\begin{equation}\label{haus}
M_\mathbf{r}^{n+1}\tau_\mathbf{r}^{-1}(M_\mathbf{r}^{-n}\mathbf{t})
\subset M_\mathbf{r}^{n+1}\tau_\mathbf{r}^{-n-1}(\mathbf{x}).
\end{equation}
Note that $\tau_\mathbf{r}$ is surjective since
$0<|r_0|=|\!\det M_\mathbf{r}|<1$. Thus there exists some
$\mathbf{t}'\in\tau_\mathbf{r}^{-1}(M_\mathbf{r}^{-n}\mathbf{t})$.
By the definition of $\tau_\mathbf{r}$, there is a
$\mathbf{v}=(0,\ldots,0,v)$ with $v\in[0,1)$ such that
$\tau_\mathbf{r}(\mathbf{t}')=M_\mathbf{r}\mathbf{t}' + \mathbf{v}$.
Now we have $M_\mathbf{r}\mathbf{t}' + \mathbf{v} =
M_\mathbf{r}^{-n}\mathbf{t}$. Using (\ref{haus}), this gives
$\mathbf{t} -
M_\mathbf{r}^{n}\mathbf{v}=M_\mathbf{r}^{n+1}\mathbf{t}'\in
M_\mathbf{r}^{n+1}\tau_\mathbf{r}^{-n-1}(\mathbf{x})$.

On the other hand, let $\mathbf{t}\in
M_\mathbf{r}^{n+1}\tau_\mathbf{r}^{-n-1}(\mathbf{x})$. Then $
M_\mathbf{r}^{n}\tau_\mathbf{r}(M_\mathbf{r}^{-n-1}\mathbf{t})\in
M_\mathbf{r}^{n}\tau_\mathbf{r}^{-n}(\mathbf{x})$. As there exists a
$\mathbf{v}=(0,\ldots,0,v)$ with $v\in[0,1)$ such that
$\tau_\mathbf{r}(M_\mathbf{r}^{-n-1}\mathbf{t})=M_\mathbf{r}^{-n}\mathbf{t}+\mathbf{v}$,
we conclude that $\mathbf{t}+M_{\mathbf{r}}^n\mathbf{v} \in
M_\mathbf{r}^n\tau_\mathbf{r}^{-n}(\mathbf{x})$.

Since
$\|M_\mathbf{r}^n\mathbf{v}\| \le \tilde\rho^n\|(0,\ldots,0,1)^t\|$ we
are done.
\end{proof}

\begin{theorem}[Basic properties of SRS tiles]\label{seteq}
Let $\mathbf{r}\in \mathrm{int}(\mathcal{D}_d)$ be reduced and
$\mathbf{x}\in\mathbb{Z}^d$. The SRS tile
$\mathcal{T}_\mathbf{r}(\mathbf{x})$ can be written as
\[
\mathcal{T}_\mathbf{r}(\mathbf{x}) = \mathop{\rm
Lim}_{n\rightarrow\infty} M_\mathbf{r}^n
\tau_\mathbf{r}^{-n}(\mathbf{x})
\]
where $\mathop{\rm Lim}$ denotes the limit w.r.t.\ the Hausdorff
metric~$\delta$.
It is a non-empty compact set that satisfies the set equation
\begin{equation} \label{eq:decomposition}
\mathcal{T}_\mathbf{r}(\mathbf{x})=\bigcup_{\mathbf{y}\in\tau_\mathbf{r}^{-1}(\mathbf{x})}M_\mathbf{r}\mathcal{T}_\mathbf{r}(\mathbf{y}).
\end{equation}
\end{theorem}

\begin{proof}
The fact that the Hausdorff limit $\mathop{\rm
Lim}_{n\rightarrow\infty} M_\mathbf{r}^n
\tau_\mathbf{r}^{-n}(\mathbf{x})$ exists and equals
$\mathcal{T}_\mathbf{r}(\mathbf{x})$ follows from
Lemma~\ref{lem:Hausdorff}. Moreover,
$\mathcal{T}_\mathbf{r}(\mathbf{x})$ is closed since Hausdorff
limits are closed by definition. As it is also bounded by
Lemma~\ref{endlich} the compactness of
$\mathcal{T}_\mathbf{r}(\mathbf{x})$ is shown.
The fact that $\mathcal{T}_\mathbf{r}(\mathbf{x})$ is non-empty follows from the surjectivity of $\tau_\mathbf{r}$.
It remains to prove
the set equation. This follows from
\begin{align*}
\mathcal{T}_\mathbf{r}(\mathbf{x}) &= \mathop{\rm
Lim}_{n\rightarrow\infty} M_\mathbf{r}^n
\tau_\mathbf{r}^{-n}(\mathbf{x}) = M_\mathbf{r}\mathop{\rm
Lim}_{n\rightarrow\infty}
\bigcup_{\mathbf{y}\in\tau_\mathbf{r}^{-1}(\mathbf{x})}
M_\mathbf{r}^{n-1} \tau_\mathbf{r}^{-n+1}(\mathbf{y}) \\
&= M_\mathbf{r}
\bigcup_{\mathbf{y}\in\tau_\mathbf{r}^{-1}(\mathbf{x})} \mathop{\rm
Lim}_{n\rightarrow\infty}M_\mathbf{r}^{n-1}
\tau_\mathbf{r}^{-n+1}(\mathbf{y})=M_{\mathbf{r}}\bigcup_{\mathbf{y}\in\tau_\mathbf{r}^{-1}(\mathbf{x})}
\mathcal{T}_\mathbf{r}(\mathbf{y}). \hfill\qedhere
\end{align*}
\end{proof}

The points in an SRS tile are characterized by the following
proposition.

\begin{proposition}\label{prop:rep}
Let $\mathbf{r} \in \mathrm{int}(\mathcal{D}_d)$ be reduced and $\mathbf{x} \in
\mathbb{Z}^d$. Then
$\mathbf{t}\in\mathcal{T}_\mathbf{r}(\mathbf{x})$ if and only if
there exist some numbers $v_{-j}\in[0,1)$, $j\in\mathbb{N}$, such
that
\begin{equation} \label{eq:txvn}
\mathbf{t} = \mathbf{x} - \sum_{j=0}^\infty
M_\mathbf{r}^j(0,\ldots,0,v_{-j})^t
\end{equation}
and
\begin{equation} \label{eq:inZd}
 M_\mathbf{r}^{-n} \Big(\mathbf{x} - \sum_{j=0}^{n-1} M_\mathbf{r}^j(0,\ldots,0,v_{-j})^t\Big) \in \mathbb{Z}^d \quad\mbox{for all }n\in\mathbb{N}.
\end{equation}
Set $\mathbf{z}_{-n} :=M_\mathbf{r}^{-n} \Big(\mathbf{x} -
\sum_{j=0}^{n-1} M_\mathbf{r}^j(0,\ldots,0,v_{-j})^t\Big)$. Then
condition (\ref{eq:inZd}) is equivalent to
$\mathbf{z}_0=\mathbf{x}$,
\begin{equation} \label{eq:tauzn}
\tau_\mathbf{r}(\mathbf{z}_{-n})=\mathbf{z}_{-n+1}
\quad\mbox{and}\quad v_{-n+1}=\{\mathbf{r}\mathbf{z}_{-n}\}
\quad\mbox{for all } n\ge 1.
\end{equation}
\end{proposition}

\begin{proof}
The equivalence of (\ref{eq:inZd}) and (\ref{eq:tauzn}) follows from
Proposition~\ref{ex:tau}. If these conditions hold and $\mathbf{t}$
is defined by (\ref{eq:txvn}), then it is clear from Remark~\ref{r:anne} that
$\mathbf{t}\in\mathcal{T}_\mathbf{r}(\mathbf{x})$.

Now let $\mathbf{t}\in\mathcal{T}_\mathbf{r}(\mathbf{x})$. We show
that we can choose the $\mathbf{z}_{-n}$ given in Remark~\ref{r:anne} such that
$\tau_\mathbf{r}(\mathbf{z}_{-n})=\mathbf{z}_{-n+1}$ for all
$n\ge1$. Let $\mathbf{z}_0=\mathbf{x}$. By (\ref{eq:decomposition}),
there is some $\mathbf{z}_{-1}\in\tau_\mathbf{r}^{-1}(\mathbf{z}_0)$
with
$M_\mathbf{r}^{-1}\mathbf{t}\in\mathcal{T}_\mathbf{r}(\mathbf{z}_{-1})$,
and inductively
$\mathbf{z}_{-n}\in\tau_\mathbf{r}^{-1}(\mathbf{z}_{-n+1})$ with
$M_\mathbf{r}^{-n}\mathbf{t}\in\mathcal{T}_\mathbf{r}(\mathbf{z}_{-n})$
for all $n\ge1$. By Lemma~\ref{endlich}, we have
$\|M_\mathbf{r}^{-n}\mathbf{t}-\mathbf{z}_{-n}\|\le R$, thus
$\lim_{n\to\infty}M_\mathbf{r}^n\mathbf{z}_{-n}=\mathbf{t}$, and
$\mathbf{t}$ is of the form (\ref{eq:txvn}) with
$v_{-n+1}=\{\mathbf{r}\mathbf{z}_{-n}\}$.
\end{proof}

It remains to show the covering property.

\begin{proposition} \label{prop:covering}
Let $\mathbf{r} \in \mathrm{int}(\mathcal{D}_d)$ be reduced.
The family of SRS
tiles
$\{\mathcal{T}_\mathbf{r}(\mathbf{x})\mid{\mathbf{x}\in\mathbb
Z^d}\}$ is a \emph{covering} of $\mathbb{R}^d$, {\it i.e.},
\[
\bigcup_{\mathbf{x}\in\mathbb{Z}^d}\mathcal{T}_\mathbf{r}(\mathbf{x}) = \mathbb{R}^d.
\]
\end{proposition}

\begin{proof}
Set $\mathcal{C} = \bigcup_{\mathbf{x}\in\mathbb{Z}^d}\mathcal{T}_\mathbf{r}(\mathbf{x})$.
By Lemma~\ref{endlich} and the non-emptiness of $\mathcal{T}_\mathbf{r}(\mathbf{x})$, the set $\mathcal{C}$ is relatively dense in $\mathbb{R}^d$.
By (\ref{eq:decomposition}), we have $M_\mathbf{r} \mathcal{C} \subseteq \mathcal{C}$.
As $M_\mathbf{r}$ is contractive, this implies that $\mathcal{C}$ is dense in $\mathbb{R}^d$.
We conclude by noticing that the SRS tiles are compact by
Theorem~\ref{seteq} and that the family of SRS tiles
$\{\mathcal{T}_\mathbf{r}(\mathbf{x}) \mid \mathbf{x}\in\mathbb Z^d\}$ is locally finite, according to Lemma~\ref{endlich}.
\end{proof}

\subsection{Around  the origin} \label{subsec:origin}
The tile associated with $\mathbf{0}$ plays a specific role.

\begin{definition}[Central SRS tile]\label{srstdefcentral}
Let $\mathbf{r}\in \mathrm{int}(\mathcal{D}_d)$ be reduced. The tile
$\mathcal{T}_\mathbf{r}(\mathbf{0})$ is called \emph{central SRS
tile} associated with~$\mathbf{r}$.
\end{definition}

Since $\tau_\mathbf{r}(\mathbf{0})=\mathbf{0}$ for every
$\mathbf{r}\in\mathrm{int}(\mathcal{D}_d)$, the origin is an element
of the central tile. However, in general it can be contained in
finitely many other tiles of the collection
$\{\mathcal{T}_\mathbf{r}(\mathbf{x}) \mid \mathbf{x}\in
\mathbb{Z}^d\}$. Whether or not $\mathbf{0}$ is contained
exclusively in the central tile plays an important role in
numeration. Indeed, for beta-numeration, $\mathbf{0}$ is contained
exclusively in the central beta-tile (see
Definition~\ref{def:fracbeta} below) if and only if the so-called
finiteness property (F) is satisfied (see
\cite{Akiyama:02,Frougny-Solomyak:92}). An analogous criterion holds
for CNS ({\it cf.}~\cite{Akiyama-Thuswaldner:00}). Now, we show that
this characterizes also the SRS with finiteness property.

\begin{definition}[Purely periodic point]
Let $\mathbf{r} \in \mathcal{D}_d$.
A point $\mathbf{z} \in \mathbb{Z}^d$ is called \emph{purely periodic} if $\tau_\mathbf{r}^p(\mathbf{z})=\mathbf{z}$ for some $p\ge 1$.
\end{definition}

\begin{theorem}\label{t63}
Let $\mathbf{r} \in \mathrm{int}(\mathcal{D}_d)$ be reduced and $\mathbf{x} \in \mathbb{Z}^d$.
Then $\mathbf{0} \in \mathcal{T}_\mathbf{r}(\mathbf{x})$ if and only if $\mathbf{x}$ is purely periodic.
There are only finitely many purely periodic points.
\end{theorem}

\begin{proof}
We first show that, if $\mathbf{x}$ is purely periodic with period
$p$, then $\mathbf{0} \in \mathcal{T}_\mathbf{r}(\mathbf{x})$. We
have $\tau_\mathbf{r}^p(\mathbf{x})=\mathbf{x}$. Therefore
$\mathbf{x} \in \tau_\mathbf{r}^{-k p}(\mathbf{x})$ for all $k \in
\mathbb{N}$, and since $M_\mathbf{r}$ is contractive we gain
\[
\mathbf{0}=\lim_{k\to\infty}M_\mathbf{r}^{k p}\mathbf{x} \in \mathcal{T}_\mathbf{r}(\mathbf{x}).
\]

To prove the other direction, let $\mathbf{0} \in
\mathcal{T}_\mathbf{r}(\mathbf{x})$, $\mathbf{x} \in \mathbb{Z}^d$.
Let $\mathbf{z}_{-n}\in\tau_\mathbf{r}^{-n}(\mathbf{x})$ be defined as
in Proposition~\ref{prop:rep}, with $\mathbf{t}=\mathbf{0}$.
Multiplying (\ref{eq:txvn}) by $M_\mathbf{r}^{-n}$, we gain
$\mathbf{0}\in \mathcal{T}_\mathbf{r}(\mathbf{z}_{-n})$ for all
$n\in\mathbb{N}$ by Proposition~\ref{prop:rep} because the expression in (\ref{eq:inZd}) is zero for each $n\in\mathbb{N}$ in this case.
By Lemma~\ref{endlich}, the set $\{\mathbf{z}_{-n} \mid n\in\mathbb{N}\}$ is finite, therefore we have $\mathbf{z}_{-n}=\mathbf{z}_{-k}$ for some $n>k\ge0$.
Since $\tau_\mathbf{r}^{n-k}(\mathbf{z}_{-n})=\mathbf{z}_{-k}$, we gain that $\mathbf{z}_{-n}$ is purely periodic, and thus the same is true for $\mathbf{x}=\tau_\mathbf{r}^n(\mathbf{z}_{-n})$.

Again, by Lemma~\ref{endlich}, it follows that only points $\mathbf{x}\in\mathbb{Z}^d$ with $\|\mathbf{x}\|\le R$ can be purely periodic.
Note that this was already proved in \cite{Akiyama-Borbeli-Brunotte-Pethoe-Thuswaldner:05}.
\end{proof}

A point $\mathbf{t}\in\mathbb{R}^d$ that satisfies $\mathbf{t} \in
\mathcal{T}_\mathbf{r}(\mathbf{x}) \setminus \bigcup_{\mathbf{y}\neq
\mathbf{x}} \mathcal{T}_\mathbf{r}(\mathbf{y})$ for some
$\mathbf{x}\in\mathbb{Z}^d$ is called an \emph{exclusive
point} of $\mathcal{T}_\mathbf{r}(\mathbf{x})$ according to the terminology introduced in \cite{Akiyama:02}
in the case of beta-tiles. Note that every exclusive point of
$\mathcal{T}_\mathbf{r}(\mathbf{x})$ is an inner point of
$\mathcal{T}_\mathbf{r}(\mathbf{x})$ because SRS tiles are compact
(Theorem~\ref{seteq}) and because the collection
$\{\mathcal{T}_\mathbf{r}(\mathbf{x})\mid \mathbf{x}\in
\mathbb{Z}^d\}$ is locally finite (Lemma~\ref{endlich}) and covers~$\mathbb{R}^d$ (Proposition~\ref{prop:covering}).
We will come back
to the notion of exclusive points in Section~\ref{sec:multiple}. The
following corollary is a consequence of
Theorem~\ref{t63}.

\begin{corollary}\label{cor:exc}
Let $\mathbf{r} \in \mathrm{int}(\mathcal{D}_d)$ be reduced.
Then $\mathbf{r} \in \mathcal{D}_d^{(0)}$ if and only if $\mathbf{0} \in \mathcal{T}_\mathbf{r}(\mathbf{0}) \setminus
\bigcup_{\mathbf{y}\neq \mathbf{0}} \mathcal{T}_\mathbf{r}(\mathbf{y})$.
\end{corollary}

\begin{proof}
Note that $\mathbf{r} \in \mathcal{D}_d^{(0)}$ if and only if each orbit ends up in $\mathbf{0}$, implying that $\mathbf{0}$ is the only purely periodic point.
\end{proof}

It immediately follows that for $\mathbf{r} \in \mathcal{D}_d^{(0)}$ the central tile $\mathcal{T}_\mathbf{r}(\mathbf{0})$ has non-empty interior.
One may ask if the interior of $\mathcal{T}_\mathbf{r}(\mathbf{x})$ is non-empty for each choice $\mathbf{r} \in \mathrm{int}(\mathcal{D}_d)$, $\mathbf{x} \in \mathbb{Z}^d$.
The answer is no, as the following example shows.

\begin{example}\label{ex:pp}
Set $\mathbf{r}=(\frac{9}{10},-\frac{11}{20})$.
Consider the points
\[
\mathbf{z}_1=(-1,-1)^t,\ \mathbf{z}_2=(-1,1)^t,\ \mathbf{z}_3=(1,2)^t,\ \mathbf{z}_4=(2,1)^t,\ \mathbf{z}_5=(1,-1)^t.
\]
It can easily be verified that
\[
\tau_\mathbf{r}:\mathbf{z}_1 \mapsto \mathbf{z}_2 \mapsto \mathbf{z}_3 \mapsto \mathbf{z}_4 \mapsto \mathbf{z}_5 \mapsto \mathbf{z}_1.
\]
Thus, each of these points is purely periodic.
Now calculate $\tau_\mathbf{r}^{-1}(\mathbf{z}_1)$:
\[
\tau_\mathbf{r}^{-1}(\mathbf{z}_1) = \Big\{(x,-1)^t \;\Big|\; x \in \mathbb{Z},\, 0 \leq \frac{9}{10} x +\frac{11}{20} -1 <1\Big\} = \big\{(1,-1)^t\big\}=\{\mathbf{z}_5\}.
\]
Similarly it can be shown that $\tau_\mathbf{r}^{-1}(\mathbf{z}_i)=\{\mathbf{z}_{i-1}\}$ for $i \in \{2,3,4,5\}$.
Hence every tile $\mathcal{T}_\mathbf{r}(\mathbf{z}_i)$, $i\in \{1,2,3,4,5\}$, consists of a single point (the point $\mathbf{0}$).
The central tile $\mathcal{T}_\mathbf{r}(\mathbf{0})$ for this parameter is the one shown in Figure~\ref{fracs} on the lower right hand side.
\end{example}

\begin{example}\label{ex:bt}
For $\mathbf{r}=(\frac{3}{4},1)$, the tiles $\mathcal{T}_\mathbf{r}(\mathbf{x})$ with $\|\mathbf{x}\|_\infty \leq 2$ are depicted in Figure~\ref{map}.
As $\mathbf{r} \in \mathcal{D}_d^{(0)}$, we will see in Corollary~\ref{cor:weaktiling} that the collection $\{\mathcal{T}_\mathbf{r}(\mathbf{x}) \mid \mathbf{x} \in \mathbb{Z}^d\}$ is a weak tiling of $\mathbb{R}^d$.
\begin{figure}[h]
\centering
\includegraphics[width=0.7\textwidth]{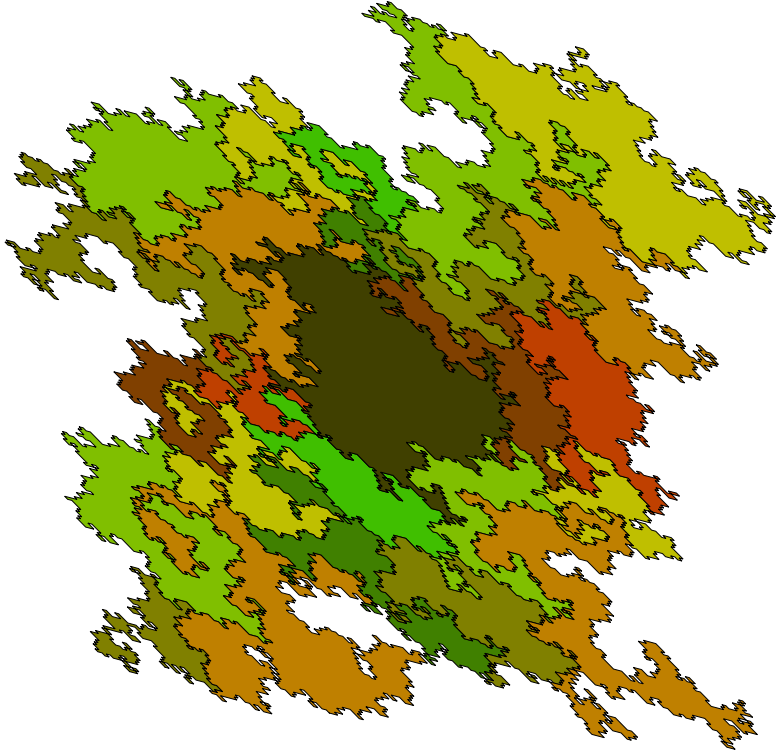}
\caption{The SRS tiles $\mathcal{T}_\mathbf{r}(\mathbf{x})$ with
$\|\mathbf{x}\|_\infty\le 2$ corresponding to
$\mathbf{r}=(\frac{3}{4},1)$.} \label{map}
\end{figure}
\end{example}

\section{Multiple tilings and tiling conditions} \label{sec:multiple}

According to Proposition~\ref{prop:covering}, the family $\{\mathcal{T}_\mathbf{r}(\mathbf{x})\mid \mathbf{x}\in
\mathbb{Z}^d\}$ of SRS tiles forms a covering of~$\mathbb{R}^d$.
Various tiling conditions concerning CNS and beta-tiles are spread  in the
literature (see \emph{e.g.}\ the references
in~\cite{Barat-Berthe-Liardet-Thuswaldner:06,Berthe-Siegel:05}). They are of a combinatorial or
dynamical nature, or they are expressed in terms of number systems.
Among these conditions, the fact that $\mathbf{0}$ is an inner
point of the central tile plays an important role,  which is related
to the finiteness property (F) introduced and discussed in
Sections~\ref{beta} (see \emph{e.g.}\
\cite{Akiyama:98} for  the case of beta-tiles).

In this section, we study tiling properties of SRS tiles. The
notions of covering and tiling we will use here are discussed   in
Section~\ref{subsec:mtiling}. In Section~\ref{subsec:approx}, some
facts on $m$-exclusive points are shown. A sufficient condition for
coverings to be in fact tilings is given  in
Section~\ref{subsec:tiling}.

\subsection{Coverings and tilings}\label{subsec:mtiling}
According to Lemma~\ref{endlich} and Proposition~\ref{prop:covering}, the family of SRS tiles $\{\mathcal{T}_\mathbf{r}(\mathbf{x}) \mid \mathbf{x}\in\mathbb{Z}^d\}$ is a \emph{covering with bounded degree}, that is, every point $\mathbf{t}\in\mathbb{R}^d$ is contained in a finite and uniformly bounded number of tiles.
Thus there exists a unique positive integer $m$ such that every point $\mathbf{t}\in\mathbb{R}^d$ is contained in at least $m$ SRS  tiles and there exists a point that is contained in exactly $m$ tiles.
Let us introduce several definitions concerning the notions of covering and tiling.

\begin{definition}[Covering and tiling; exclusive and inner point]\label{def:tiling}
Let $\mathcal{K}$ be a locally finite collection of compact subsets covering $\mathbb{R}^d$.

\begin{itemize}
\item
The \emph{covering degree} with respect to $\mathcal{K}$ of a point $\mathbf{t}\in\mathbb{R}^d$ is given by $\deg(\mathcal{K},\mathbf{t}):=\#\{K\in\mathcal{K} \mid \mathbf{t}\in K\}$.
\item
The \emph{covering degree} of $\mathcal{K}$ is given by $\deg(\mathcal{K}) := \min \{\deg(\mathcal{K},\mathbf{t}) \mid \mathbf{t}\in\mathbb{R}^d\}$.
\item
The collection $\mathcal{K}$ is a \emph{weak $m$-tiling} if $\deg(\mathcal{K})=m$ and $\bigcap_{i=1}^{m+1}\mathrm{int}(K_i)=\emptyset$ for every choice of $m+1$ pairwise different $K_1,\ldots,K_{m+1}\in\mathcal{K}$.
A~weak $1$-tiling is also called \emph{weak tiling}.
\item
A point $\mathbf{t}\in\mathbb{R}^d$ is \emph{$m$-exclusive} with respect to $\mathcal{K}$ if $\deg(\mathcal{K},\mathbf{t}) = \deg(\mathcal{K}) = m$.
\item
A point $\mathbf{t}\in\mathbb{R}^d$ is an \emph{inner point} of the collection $\mathcal{K}$ if  $\mathbf{t}\not\in\bigcup_{K\in\mathcal{K}}\partial K$.
\end{itemize}
\end{definition}

In particular, the collection $\mathcal{K}$ is a weak tiling if  each inner point of $\mathcal{K}$ belongs to exactly one element of~$\mathcal{K}$.
Moreover, the definition of $1$-exclusive points with respect to $\{\mathcal{T}_{\mathbf{r}}(\mathbf{x}) \mid {\mathbf{x}\in\mathbb{Z}^d}\}$ recovers the notion of exclusive points introduced in Section~\ref{subsec:origin}.

Let us recall that a tiling by translation is often defined as a    collection of tiles having finitely many tiles up to translation, with a tile being assumed to be the closure of its interior.
We study weak tilings in the sense of Definition~\ref{def:tiling} because of the following reasons:
\begin{itemize}
\item
There exist choices of $\mathbf{r}\in\mathrm{int}(\mathcal{D}_d)$,
$\mathbf{x}\in\mathbb{Z}^d$, such that the tile $\mathcal{T}_\mathbf{r}(\mathbf{x})$ is not the closure of its interior, see Example~\ref{ex:pp}.
\item
There exist parameters $\mathbf{r}\in\mathrm{int}(\mathcal{D}_d)$ such that the family $\{\mathcal{T}_{\mathbf{r}}(\mathbf{x}) \mid {\mathbf{x}\in\mathbb{Z}^d}\}$ is not a collection of finitely many tiles up to translation, \emph{e.g.}\ $\mathbf{r}=(-2/3)$, see Corollary~\ref{cor:515}.
We conjecture that this holds for every $\mathbf{r}$ related to a non-monic CNS (see Section~\ref{cns}) or a non-unit Pisot number (see Section~\ref{beta}).
\item
We are not able to show that the boundaries of the tiles
$\mathcal{T}_{\mathbf{r}}(\mathbf{x})$ have zero $d$-dimensional
Lebesgue measure. If $\mathbf{r}$ is related to a monic CNS or a
unit Pisot number, the fact that the boundary of each tile has zero
measure is a direct consequence of the self-affine structure of the
boundary of the tiles (\emph{cf.}\
\cite{Akiyama:02,Kalle-Steiner:09,Lagarias-Wang:96a}). For other
parameters $\mathbf{r}$, we have no appropriate description of the
boundary, and we cannot evaluate its measure.
\end{itemize}

We are now going to prove that, for a large class of parameters~$\mathbf{r}$, the collection $\{\mathcal{T}_{\mathbf{r}}(\mathbf{x}) \mid {\mathbf{x}\in\mathbb{Z}^d}\}$ is a weak $m$-tiling, by showing that the set of $m$-exclusive points is dense in~$\mathbb{R}^d$.

\subsection{$m$-exclusive points} \label{subsec:approx}

Let us first prove that $m$-exclusive points are always inner points of $\{\mathcal{T}_{\mathbf{r}}(\mathbf{x}) \mid {\mathbf{x}\in\mathbb{Z}^d}\}$.

\begin{lemma}\label{l:32}
Let $\mathbf{r} \in \mathrm{int}(\mathcal{D}_d)$ be reduced and let $m$ be the covering degree of $\{\mathcal{T}_\mathbf{r}(\mathbf{x}) \mid \mathbf{x}\in\mathbb{Z}^d\}$.
Then there exists an $m$-exclusive point.
Every $m$-exclusive point $\mathbf{t}\in \bigcap_{i=1}^m\mathcal{T}_\mathbf{r}(\mathbf{x}_i)$ satisfies $\mathbf{t} \in \bigcap_{i=1}^m\mathrm{int}(\mathcal{T}_\mathbf{r}(\mathbf{x}_i))$.
In particular, $m$-exclusive points are inner points of $\{\mathcal{T}_\mathbf{r}(\mathbf{x}) \mid \mathbf{x}\in\mathbb{Z}^d\}$.
\end{lemma}

\begin{proof}
The existence of an $m$-exclusive point is a direct consequence of the definition of the covering degree.

Assume that $\mathbf{t}$ is an $m$-exclusive point.
Let $\mathcal{T}_\mathbf{r}(\mathbf{x}_1),\ldots,\mathcal{T}_\mathbf{r}(\mathbf{x}_m)$ be the $m$ tiles it belongs to.
Assume that there exists a sequence of points $(\mathbf{t}_n)_{n\in\mathbb{N}}$ with values in $\mathbb{R}^d$ such that $\lim_{n\to\infty}\mathbf{t}_n = \mathbf{t}$, and a sequence of points $(\mathbf{z}_n)_{n\in\mathbb{N}}$ with values in $\mathbb{Z}^d$ such that
$\mathbf{t}_n\in\mathcal{T}_\mathbf{r}(\mathbf{z}_n)$ for all $n\in\mathbb{N}$, with $\mathbf{z}_n$ distinct from all $\mathbf{x}_i$'s.
Since the collection is locally finite, an infinite subsequence of $\mathbf{z}_n$'s is constant, say equal to $\mathbf{z}$.
Since the tiles are compact, this implies that $\mathbf{t} \in \mathcal{T}_\mathbf{r}(\mathbf{z})$, which contradicts the $m$-exclusivity.
Thus there exists a neighbourhood $U$ of $\mathbf{t}$ with $U \cap \mathcal{T}_\mathbf{r}(\mathbf{z}) = \emptyset$ for each $\mathbf{z} \in \mathbb{Z}^d \setminus \{\mathbf{x}_1,\ldots,\mathbf{x}_m\}$.
By the $m$-covering property, we know that each point belongs to at
least $m$ tiles, which implies that $U \subseteq \bigcap_{i=1}^m\mathcal{T}_\mathbf{r}(\mathbf{x}_i)$.
\end{proof}

Let us now prove that an $m$-exclusive point $\mathbf{t} \in \bigcap_{i=1}^m\mathcal{T}_\mathbf{r}(\mathbf{x}_i)$ is somehow characterized by any sequence of approximations $M_\mathbf{r}^n\mathbf{z}_{-n}$ defined in Proposition~\ref{prop:rep}.
Note that $\mathbf{t} = \mathop{\rm Lim}_{n\to\infty} M_\mathbf{r}^n\mathcal{T}_\mathbf{r}(\mathbf{z}_{-n})$, and that $M_\mathbf{r}^n\mathcal{T}_\mathbf{r}(\mathbf{z}_{-n})$ is a tile in the $n$-fold subdivision $\mathcal{T}_\mathbf{r}(\mathbf{x}_i) = \bigcup_{\mathbf{z}\in\tau_\mathbf{r}^{-n}(\mathbf{x}_i)} M_\mathbf{r}^n \mathcal{T}_\mathbf{r}(\mathbf{z})$ for some~$\mathbf{x}_i$, which is given by applying (\ref{eq:decomposition}) $n$ times.
The following proposition states that a point $\mathbf{t}$ is $m$-exclusive provided that, for some $n\in\mathbb{N}$, each $M_\mathbf{r}^n \mathcal{T}_\mathbf{r}(\mathbf{z}_{-n}+\mathbf{y})$ with $\|\mathbf{y}\|\le 2R$ occurs in the subdivision of some tile $\mathcal{T}_\mathbf{r}(\mathbf{x}_i)$, $1\le i\le m$.

\begin{proposition}\label{prop:ex}
Let $m$ be the covering degree of $\{\mathcal{T}_\mathbf{r}(\mathbf{x}) \mid \mathbf{x}\in\mathbb{Z}^d\}$, $\mathbf{t}\in\mathbb{R}^d$, $(\mathbf{z}_{-n})_{n\in\mathbb{N}}$ as in Proposition~\ref{prop:rep}, and $R$ defined by (\ref{eq:R}).
Then $\mathbf{t}$ is $m$-exclusive and contained in the intersection $\bigcap_{i=1}^m \mathcal{T}_\mathbf{r}(\mathbf{x}_i)$ if and only if there exists some $n\in\mathbb{N}$ such that
\begin{equation} \label{eq:2R}
\tau_\mathbf{r}^n(\mathbf{z}_{-n}+\mathbf{y}) \in \{\mathbf{x}_1,\ldots,\mathbf{x}_m\} \quad\mbox{for all}\ \mathbf{y}
\in \mathbb{Z}^d\ \mbox{with}\ \|\mathbf{y}\|\le 2R.
\end{equation}

If, for some $\mathbf{z}\in\mathbb{Z}^d$, $n\in\mathbb{N}$,
\begin{equation}\label{f:inclusion}
\#\{\tau_\mathbf{r}^n(\mathbf{z}+\mathbf{y}) \mid \mathbf{y}
\in \mathbb{Z}^d,\,\|\mathbf{y}\|\le R\big\} = m,
\end{equation}
then $M_\mathbf{r}^n\mathbf{z}$ is an $m$-exclusive point.
\end{proposition}

\begin{proof}
Assume that $\mathbf{t}$ is an $m$-exclusive point. By
Lemma~\ref{l:32}, there exists some $\varepsilon>0$ such that every
point $\mathbf{t}'\in\mathbb{R}^d$ satisfying
$\|\mathbf{t}'-\mathbf{t}\|<\varepsilon$ lies only in the tiles
$\mathcal{T}_{\mathbf{r}}(\mathbf{x}_1),\,\ldots,\,\mathcal{T}_{\mathbf{r}}(\mathbf{x}_m)$.
Let $n\in\mathbb{N}$ satisfy $4\tilde\rho^n R<\varepsilon$. Since
$M_\mathbf{r}^{-n}\mathbf{t} \in
\mathcal{T}_\mathbf{r}(\mathbf{z}_{-n})$, we have
$\|M_\mathbf{r}^{-n}\mathbf{t}-\mathbf{z}_{-n}\|\le R$ by
Lemma~\ref{endlich}, thus
$\|\mathbf{t}-M_\mathbf{r}^n(\mathbf{z}_{-n}+\mathbf{y})\|\le
3\tilde\rho^n R$ if $\|\mathbf{y}\|\le 2R$. By Theorem~\ref{seteq},
there exists a point $\mathbf{t}' \in
\mathcal{T}_\mathbf{r}\big(\tau_\mathbf{r}^n(\mathbf{z}_{-n}+\mathbf{y})\big)$
with $\|\mathbf{t}'-M_\mathbf{r}^n(\mathbf{z}_{-n}+\mathbf{y})\| \le
\tilde\rho^n R$. Since $\|\mathbf{t}'-\mathbf{t}\|\le 4\tilde\rho^n
R<\varepsilon$, we obtain
$\tau_\mathbf{r}^n(\mathbf{z}_{-n}+\mathbf{y}) \in
\{\mathbf{x}_1,\ldots,\mathbf{x}_m\}$.

For the other direction, assume that there exists some $n\in\mathbb{N}$ such that (\ref{eq:2R}) holds.
We have to show that $\mathbf{t} \in \mathcal{T}_{\mathbf{r}}(\mathbf{z}_0')$ implies $\mathbf{z}_0' \in \{\mathbf{x}_1,\ldots,\mathbf{x}_m\}$.
Let $(\mathbf{z}_{-n}')_{n\in\mathbb{N}}$ be as in Proposition~\ref{prop:rep}.
Since $\|M_\mathbf{r}^{-n}\mathbf{t}-\mathbf{z}_{-n}'\|\le R$ and $\|M_\mathbf{r}^{-n}\mathbf{t}-\mathbf{z}_{-n}\|\le R$, (\ref{eq:2R}) implies $\tau_\mathbf{r}^n(\mathbf{z}_{-n}') \in \{\mathbf{x}_1,\ldots,\mathbf{x}_m\}$.
Since $\mathbf{z}_0'=\tau_\mathbf{r}^n(\mathbf{z}_{-n}')$, the point $\mathbf{t}$ is $m$-exclusive.

For the second statement, let $\{\tau_\mathbf{r}^n(\mathbf{z}+\mathbf{y}) \mid \mathbf{y}
\in \mathbb{Z}^d,\,\|\mathbf{y}\|\le R\big\} = \{\mathbf{x}_1,\ldots,\mathbf{x}_m\}$, $M_\mathbf{r}^n\mathbf{z} \in \mathcal{T}_{\mathbf{r}}(\mathbf{z}_0')$, and $(\mathbf{z}_{-n}')_{n\in\mathbb{N}}$ be as in Proposition~\ref{prop:rep}, with $\mathbf{t}=M_\mathbf{r}^n\mathbf{z}$.
Then we have $\|\mathbf{z}_{-n}'-\mathbf{z}\|\le R$, thus $\mathbf{z}_0' \in \{\mathbf{x}_1,\ldots,\mathbf{x}_m\}$ and $M_\mathbf{r}^n\mathbf{z}$ is $m$-exclusive.
\end{proof}

Note that Proposition~\ref{prop:ex} provides an easy way to show that a point is $m$-exclusive.
However, it does not provide a finite method to prove that a point is not $m$-exclusive.

The following corollary of Proposition~\ref{prop:ex} permits us to obtain an $m$-exclusive point from another one, by finding a translation preserving the local configuration of tiles up to the $n$-th level, with $n$ such that (\ref{f:inclusion}) holds.

\begin{corollary}\label{cor:ex}
Let $m$ be the covering degree of $\{\mathcal{T}_\mathbf{r}(\mathbf{x}) \mid \mathbf{x}\in\mathbb{Z}^d\}$, and assume that $\mathbf{z}\in\mathbb{Z}^d$, $n\in\mathbb{N}$ satisfy (\ref{f:inclusion}).
Let $\mathbf{a}\in\mathbb{Z}^d$.
If there exists some $\mathbf{b}\in\mathbb{Z}^d$ such that
\[
\tau_\mathbf{r}^n(\mathbf{z}+\mathbf{a}+\mathbf{y}) = \tau_\mathbf{r}^n(\mathbf{z}+\mathbf{y})+\mathbf{b} \quad\mbox{for all } \mathbf{y}\in\mathbb{Z}^d \mbox{ with }\|\mathbf{y}\|\le R,
\]
then $M_\mathbf{r}^n(\mathbf{z}+\mathbf{a})$ is an $m$-exclusive point.
\end{corollary}

A simple way to obtain an infinite number of $m$-exclusive points from one $m$-exclusive point $\mathbf{t}\ne\mathbf{0}$ is provided by the following lemma.

\begin{lemma} \label{lem:Mexclusive}
If $\mathbf{t}$ is an $m$-exclusive point, then $M_\mathbf{r}\mathbf{t}$ is an $m$-exclusive point.
\end{lemma}

\begin{proof}
If $M_\mathbf{r}\mathbf{t}\in\mathcal{T}_\mathbf{r}(\mathbf{x}_i)$, then there exists some $\mathbf{y}_i\in\tau_\mathbf{r}^{-1}(\mathbf{x}_i)$ such that $\mathbf{t}\in\mathcal{T}_\mathbf{r}(\mathbf{y}_i)$, and all $\mathbf{y}_i$ are mutually different if the $\mathbf{x}_i$ are.
Therefore the number of tiles to which a point $M_\mathbf{r}\mathbf{t}$ belongs cannot be larger than that for $\mathbf{t}$.
\end{proof}

\subsection{Weak $m$-tilings} \label{subsec:tiling}
In what follows we will establish our tiling results. In order to
prove these results we first show that if $\mathbf{r}$ has certain
algebraic properties then there exists a set that is relatively
dense in $\mathbb{R}^d$ containing only vectors that stabilize the
configuration of tiles. Together with Corollary~\ref{cor:ex} and
Lemma~\ref{lem:Mexclusive}, this will prove that
$\{\mathcal{T}_{\mathbf{r}}(\mathbf{x}) \mid
\mathbf{x}\in\mathbb{Z}^d\}$ forms a weak $m$-tiling. If
$\mathbf{r}$ also satisfies the finiteness property, it even forms a
weak tiling.

\begin{theorem} \label{theo:mdense}
Let $\mathbf{r}=(r_0,\ldots,r_{d-1}) \in
\mathrm{int}(\mathcal{D}_d)$ with $r_0 \ne 0$, let $m$ be the covering degree of
$\{\mathcal{T}_{\mathbf{r}}(\mathbf{x}) \mid
\mathbf{x}\in\mathbb{Z}^d\}$, and assume that $\mathbf{r}$ satisfies
one of the following conditions:
\begin{itemize}
\item
$\mathbf{r}\in\mathbb{Q}^d$,
\item
$(x-\beta)(x^d+r_{d-1}x^{d-1}+\cdots+r_1x+r_0)\in\mathbb{Z}[x]$ for some $\beta>1$,
\item
$r_0,\ldots,r_{d-1}$ are algebraically independent over $\mathbb{Q}$.
\end{itemize}
Then the set of $m$-exclusive points is dense in $\mathbb{R}^d$, and $\{\mathcal{T}_{\mathbf{r}}(\mathbf{x}) \mid \mathbf{x}\in\mathbb{Z}^d\}$ is a weak $m$-tiling.
\end{theorem}

\begin{proof}
By the definition of the covering degree
(Definition~\ref{def:tiling}), there exists an $m$-exclusive point
$\mathbf{t}$ with respect to $\{\mathcal{T}_{\mathbf{r}}(\mathbf{x})
\mid \mathbf{x}\in\mathbb{Z}^d\}$. Thus the first part of
Proposition~\ref{prop:ex} implies that there exist
$\mathbf{z}\in\mathbb{Z}^d$ and $n\in\mathbb{N}$ satisfying
(\ref{f:inclusion}). If we find a relatively dense set $\Lambda$ of
vectors $\mathbf{a} \in \mathbb{Z}^d$ satisfying the conditions of
Corollary~\ref{cor:ex}, then the set
$\{M_\mathbf{r}^n(\mathbf{z}+\mathbf{a}) \mid \mathbf{a}\in
\Lambda\}$ forms a set of $m$-exclusive inner points which is
relatively dense in $\mathbb{R}^d$. Applying
Lemma~\ref{lem:Mexclusive} to the elements of this set yields that
the set of $m$-exclusive points is dense in $\mathbb{R}^d$. In view
of Definition~\ref{def:tiling}, this already proves that
$\{\mathcal{T}_\mathbf{r}(\mathbf{x}) \mid
\mathbf{x}\in\mathbb{Z}^d\}$ forms a weak $m$-tiling.

It remains to find a relatively dense set $\Lambda$ of vectors $\mathbf{a} \in \mathbb{Z}^d$ satisfying the conditions of Corollary~\ref{cor:ex}. This is done separately for
each of the three classes of parameters given in the statement of
the theorem.

\medskip
\noindent {\bf Case 1: $\mathbf{r} \in \mathbb{Q}^d$}. Let
$\mathbf{z}\in\mathbb{Z}^d$ and $n\in\mathbb{N}$ satisfying
(\ref{f:inclusion}) be given as above. Choose $q\in\mathbb{N}$ in a
way that $\mathbf{r} \in q^{-1}{\mathbb{Z}}^d$. Then we have for
every $\mathbf{x},\mathbf{a}\in\mathbb{Z}^d$, $k\ge1$,
\[
\tau_\mathbf{r}(\mathbf{x}+q^k\mathbf{a}) = M_\mathbf{r}
(\mathbf{x}+q^k\mathbf{a}) +
\big(0,\ldots,0,\big\{\mathbf{r}(\mathbf{x}+q^k\mathbf{a})\big\}\big)
= \tau_\mathbf{r}(\mathbf{x}) + q^{k-1}\mathbf{a}'
\]
for some $\mathbf{a}'\in\mathbb{Z}^d$ which does not depend on
$\mathbf{x}$. Iterating this, we get that
$\tau_\mathbf{r}^n(\mathbf{x}+q^n\mathbf{a}) =
\tau_\mathbf{r}^n(\mathbf{x}) + \mathbf{b}$ for some
$\mathbf{b}\in\mathbb{Z}^d$ which does not depend on $\mathbf{x}$.
In particular, this implies that
\[
\tau_\mathbf{r}^n(\mathbf{x}+q^n\mathbf{a} + \mathbf{y}) =
\tau_\mathbf{r}^n(\mathbf{x}+\mathbf{y}) + \mathbf{b} \quad\hbox{for
all $\mathbf{y} \in \mathbb{Z}^d$ with  $\|\mathbf{y}\| \le R$}.
\]
Thus each element of the set $\Lambda:=\{q^n\mathbf{a} \mid
\mathbf{a}\in \mathbb{Z}^d\}$ satisfies the conditions of
Corollary~\ref{cor:ex}. As $\Lambda$ is relatively dense in~$\mathbb{R}^d$, we are done in this case.

\medskip
\noindent {\bf Case 2:
$A(x)=(x-\beta)(x^d+r_{d-1}x^{d-1}+\cdots+r_1x+r_0)\in\mathbb{Z}[x]$
for some $\beta>1$}. Let $\mathbf{z}\in\mathbb{Z}^d$ and
$n\in\mathbb{N}$ satisfying (\ref{f:inclusion}) be given as above.
Since $\mathbf{r} \in \mathrm{int}(\mathcal{D}_d)$, all roots of
$x^d+r_{d-1}x^{d-1}+\cdots+r_1x+r_0$ have modulus less than~$1$.
Therefore, $A(x)$ is irreducible. (Indeed, $\beta$ is a Pisot
number.) Let
\[
\varepsilon := \min_{\|\mathbf{y}\|\le R,\, 0\le k<n} \beta^{-k} \big(1-\big\{\mathbf{r}\tau_\mathbf{r}^k(\mathbf{z}+\mathbf{y})\big\}\big) >0.
\]
{}From the definition of $\tau_\mathbf{r}$ we know that for every $\mathbf{x}$, $\mathbf{a} \in \mathbb{Z}^d$, one has $\tau_\mathbf{r}(\mathbf{x}+\mathbf{a}) = \tau_\mathbf{r}(\mathbf{x}) + \tau_\mathbf{r}(\mathbf{a})$ if and only if $\{\mathbf{r}\mathbf{x}+\mathbf{r}\mathbf{a}\} = \{\mathbf{r}\mathbf{x}\}+\{\mathbf{r}\mathbf{a}\}$.
In Proposition~\ref{prop:betanumformula}, we will see that $\{\mathbf{r}\tau_\mathbf{r}^k(\mathbf{a})\} = T_\beta^k(\{\mathbf{r}\mathbf{a}\})$, where $T_\beta$ is the $\beta$-transformation defined in Section~\ref{Tbeta}.
If we choose $\mathbf{a}\in\mathbb{Z}^d$ such that $\{\mathbf{r}\mathbf{a}\}<\varepsilon$, then we get $\{\mathbf{r}\tau_\mathbf{r}^k(\mathbf{a})\} = \beta^k\{\mathbf{r}\mathbf{a}\}$ for $0\le k<n$, and
\begin{equation}\label{e:a}
\tau_\mathbf{r}^k(\mathbf{z}+\mathbf{a}+\mathbf{y}) =  \tau_\mathbf{r}^k(\mathbf{z}+\mathbf{y}) + \tau_\mathbf{r}^k(\mathbf{a}) \quad\mbox{for all }\|\mathbf{y}\|\le R,\, 0\le k\le n.
\end{equation}
Thus each element of the set
$\Lambda:=\{\mathbf{a}\in\mathbb{Z}^d \mid \{\mathbf{r}\mathbf{a}\}<\varepsilon\}$
satisfies the conditions of Corollary~\ref{cor:ex}. Since $A(x)$ is
irreducible, the coordinates of $\mathbf{r}$ are linearly independent
over~$\mathbb{Q}$. Thus Kronecker's theorem yields that
$\Lambda$ is relatively dense in $\mathbb{R}^d$ and we are done
also in this case.

\medskip
\noindent {\bf Case 3: $r_0,\ldots,r_{d-1}$ are algebraically
independent}. Let $\mathbf{z}\in\mathbb{Z}^d$ and $n\in\mathbb{N}$
satisfying (\ref{f:inclusion}) be given as above. Set
\[
\varepsilon_k := \min_{\|\mathbf{y}\|\le R}
\big(1-\big\{\mathbf{r}\tau_\mathbf{r}^k(\mathbf{z}+\mathbf{y})\big\}\big)
\quad\mbox{for }0\le k<n.
\]
We have $\varepsilon_k>0$. Similarly to Case~2, each element of the
set
\[
\Lambda := \{\mathbf{a} \in \mathbb{Z}^d \mid
\{\mathbf{r}\tau_\mathbf{r}^k(\mathbf{a})\}<\varepsilon_k \hbox{ for
} 0\le k<n\}
\]
satisfies (\ref{e:a}). Thus each $\mathbf{a}\in\Lambda$ satisfies
the conditions of Corollary~\ref{cor:ex}. It remains to prove that
$\Lambda$ is relatively dense in $\mathbb{R}^d$. To this matter we
need the following notations. Let $\lfloor x\rceil$ denote the
nearest integer to $x\in\mathbb{R}$ (with $\lfloor z+1/2\rceil=z$
for $z\in\mathbb{Z}$), let $\{x\}_c=x-\lfloor x\rceil$ be the
``centralized fractional part'' and set
\[
\widetilde\tau_\mathbf{r}(\mathbf{x}) = M_\mathbf{r}\mathbf{x} +
(0,\ldots,0,\{\mathbf{r}\mathbf{x}\}_c).
\]
Consider the partition $\Lambda = \bigcup_{\eta \in \{0,1\}^n}
\Lambda_\eta$ with
\[
\Lambda_\eta := \{\mathbf{a} \in \mathbb{Z}^d \mid
\{\mathbf{r}\tau_\mathbf{r}^k(\mathbf{a})\}\in
\eta_k\varepsilon_k/2+[0,\varepsilon_k/2) \hbox{ for } 0\le k<n\}
\qquad(\eta \in \{0,1\}^n)
\]
(here $\eta=(\eta_0,\ldots,\eta_{n-1})$). We will prove the
following claim.

\noindent{\it Claim.} Let $\mathbf e_j$, $1\le j\le d$, be the canonical unit
vectors. For each $\eta \in \{0,1\}^n$ and $j\in\{1,\ldots,d\}$,
there exists $z_{\eta, j} \in \mathbb{Z}$ with
\begin{equation}\label{2ndclaim}
\big\{\mathbf{r}\widetilde\tau_\mathbf{r}^k(z_{\eta, j}\mathbf{e}_j)\}_c
\in (-1)^{\eta_k}[0,\varepsilon_k/2)  \quad\mbox{for all}\ 0\le k<n.
\end{equation}

Before we prove this claim we show that it implies the relative
denseness of $\Lambda$ in $\mathbb{R}^d$. Indeed, let $\eta \in
\{0,1\}^n$ and $j\in\{1,\ldots,d\}$ be arbitrary. Then the claim
yields that for each $\mathbf{a}\in\Lambda_\eta$ we have
\[
\big\{\mathbf{r}\tau_\mathbf{r}^k(\mathbf{a}+z_{\eta, j}\mathbf{e}_j)\big\} = \big\{\mathbf{r}\big(\tau_\mathbf{r}^k(\mathbf{a})+\widetilde\tau_\mathbf{r}^k(z_{\eta, j}\mathbf{e}_j)\big)\big\} = \big\{\mathbf{r}\tau_\mathbf{r}^k(\mathbf{a})\}+\{\mathbf{r}\widetilde\tau_\mathbf{r}^k(z_{\eta, j}\mathbf{e}_j)\big\}_c\in [0,\varepsilon_k)
\]
and thus
$\mathbf{a}+z_{\eta, j}\mathbf{e}_j \in \Lambda$. Moreover, by
analogous reasoning we see that the claim implies the following
``dual'' result: for each $\eta \in \{0,1\}^n$ set
$\eta':=(1,\ldots,1) - \eta\in \{0,1\}^n$. Then for each
$j\in\{1,\ldots,d\}$ and each $\mathbf{a} \in \Lambda_\eta$ we have
$\mathbf{a} - z_{\eta', j}\mathbf{e}_j \in \Lambda$.

Summing up, the claim implies that from each $\mathbf{a} \in
\Lambda$ there exist other elements of $\Lambda$ in uniformly
bounded distance in all positive and negative coordinate directions.
This proves that $\Lambda$ is relatively dense in $\mathbb{R}^d$.
Thus it remains to prove the above claim.

To prove this claim let $\eta=(\eta_0,\ldots,\eta_{n-1}) \in
\{0,1\}^n$ and $j\in\{1,\ldots, d\}$ be arbitrary but fixed. We need
to find an integer $z_{\eta, j}$ satisfying \eqref{2ndclaim}. First
observe that for $\mathbf{b}\in\mathbb{Z}^d$
\begin{equation} \label{eq:sfrac}
\widetilde\tau_\mathbf{r}^k(\mathbf{b}) = M_\mathbf{r}^k\mathbf{b} +
\sum_{j=0}^{k-1} M_\mathbf{r}^{k-j-1}
\big(0,\ldots,0,\{\mathbf{r}\widetilde\tau_{\mathbf{r}}^j(\mathbf{b})\}_c\big)^t.
\end{equation}
Let $\gamma_k=\{\mathbf{r}M_\mathbf{r}^k\mathbf e_d\}_c$, $0\le k<n$.
If the arguments of all centralized fractional parts are small, then
multiplying (\ref{eq:sfrac}) by $\mathbf{r}$ and applying
$\{\cdot\}_c$ yields
\begin{equation}\label{eq:sfrac2}
\{\mathbf{r}\widetilde\tau_{\mathbf{r}}^k(\mathbf{b})\}_c =
\{\mathbf{r}M_\mathbf{r}^k\mathbf{b}\}_c + \sum_{j=0}^{k-1}
\{\mathbf{r}\widetilde\tau_{\mathbf{r}}^j(\mathbf{b})\}_c\,
\gamma_{k-j-1}.
\end{equation}
Inserting (\ref{eq:sfrac2}) iteratively in itself we gain
\begin{align*}
\{\mathbf{r}\widetilde\tau_{\mathbf{r}}^k(\mathbf{b})\}_c & =   \{\mathbf{r}M_\mathbf{r}^k\mathbf{b}\}_c
+ \sum_{k=j_0>j_1>\cdots>j_\ell\ge 0,\,\ell\ge 1} \{\mathbf{r}M_\mathbf{r}^{j_\ell}\mathbf{b}\}_c \prod_{i=0}^{\ell-1} \gamma_{j_i-j_{i+1}-1} \\
& = \sum_{h=0}^k \{\mathbf{r}M_\mathbf{r}^h\mathbf{b}\}_c
\sum_{k=j_0>\cdots>j_\ell=h,\,\ell\ge0} \prod_{i=0}^{\ell-1}
\gamma_{j_i-j_{i+1}-1}.
\end{align*}
Setting
\[
c_k := \sum_{k=j_0>\cdots>j_\ell=0,\,\ell\ge0} \prod_{i=0}^{\ell-1}
\gamma_{j_i-j_{i+1}-1}
\]
we get
\[
\{\mathbf{r}\widetilde\tau_{\mathbf{r}}^k(\mathbf{b})\}_c =
\sum_{h=0}^k c_{k-h}\,\{\mathbf{r}M_\mathbf{r}^h\mathbf{b}\}_c.
\]
Now we inductively choose intervals $I_k$, $0\le k<n$, satisfying
\[
I_k+\sum_{h=0}^{k-1}c_{k-h}I_h \subseteq (-1)^{\eta_k}[0,\varepsilon_k/2).
\]
W.l.o.g., we can choose the intervals $I_k$ sufficiently small such that (\ref{eq:sfrac2}) holds provided that $\{\mathbf{r}M_\mathbf{r}^k\mathbf{b}\}_c \in I_k$ for $0\le k<n$.

Since $r_0,\ldots,r_{d-1}$ are algebraically independent, the
numbers $\{\mathbf{r}M_\mathbf{r}^k\mathbf e_j\}_c$, for $0\le k<n$,
are linearly independent over $\mathbb Q$. Thus Kronecker's theorem
yields the existence of an integer $z_{\eta, j}$ satisfying
$\{z_{\eta, j}\mathbf{r}M_\mathbf{r}^k\mathbf e_j\}_c \in I_k$ for
$0\le k<n$, hence, $z_{\eta, j}$ satisfies \eqref{2ndclaim} and we
are done.
\end{proof}

Since Corollary~\ref{cor:exc} implies that for each SRS with
finiteness property the origin is an exclusive point of
$\mathcal{T}_\mathbf{r}(\mathbf{0})$ we gain the following tiling
property.

\begin{corollary} \label{cor:weaktiling}
If $\mathbf r\in\mathcal D_d^{(0)}$ is reduced and satisfies one of the conditions
of Theorem~\ref{theo:mdense}, then
$\{\mathcal{T}_{\mathbf{r}}(\mathbf{x}) \mid
\mathbf{x}\in\mathbb{Z}^d\}$ is a weak tiling.
\end{corollary}

Let us stress the fact that we have no general algorithmic criterion
to check Proposition~\ref{prop:ex}. This is mainly due to the fact
that we have no IFS describing the boundary of an SRS tile.
Nonetheless, Theorem~\ref{theo:mdense} implies a tiling criterion.

\begin{corollary}
If $\mathbf{r} \in \mathrm{int}(\mathcal D_d)$ is reduced and satisfies one of the conditions
of Theorem~\ref{theo:mdense}, then the collection
$\{\mathcal{T}_{\mathbf{r}}(\mathbf{x}) \mid
\mathbf{x}\in\mathbb{Z}^d\}$ is a weak tiling if and only if it has at least one
exclusive point.
\end{corollary}

We conclude this section with a result that treats the case $d=1$ in
a very complete way. (We identify one dimensional vectors with
scalars here.)

\begin{theorem}\label{d1}
Let $r \in \mathrm{int}(\mathcal{D}_1)$ be reduced, \emph{i.e.},
$0<|r|<1$. Then $\{\mathcal{T}_{r}(x) \mid x\in\mathbb{Z}\}$ is a
tiling of~$\mathbb{R}$ by (possibly degenerate) intervals. Here,
\emph{tiling} has to be understood in the usual sense, \emph{i.e.},
\[
\bigcup_{x\in \mathbb{Z}} \mathcal{T}_{r}(x) = \mathbb {R} \quad
\hbox{with} \quad \#\big(\mathcal{T}_{r}(x) \cap \mathcal{T}_{r}(x')\big) \le 1 \quad\mbox{for}\ x,x'\in\mathbb{Z},\ x \neq x'.
\]
\end{theorem}

\begin{proof}
Let first $r>0$ and $x_0,y_0\in \mathbb{Z}$. Then by the definition
of $\tau_r$ we easily see that $x_0
> y_0$ implies that $-x_1 > -y_1$ for each $x_1 \in \tau_r^{-1}(x_0)$,
$y_1 \in \tau_r^{-1}(y_0)$. Thus, by induction on $n$ we have
\begin{equation}\label{1dim}
x_0 > y_0 \quad\hbox{implies that}\quad (-1)^nx_n > (-1)^ny_n
\quad\hbox{for each}\quad x_n \in \tau_r^{-n}(x_0),\, y_n \in
\tau_r^{-n}(y_0).
\end{equation}
Observe that $M_r=-r$ holds for the companion matrix in this case. Renormalizing \eqref{1dim} we now get the following assertion. Suppose that  $x_0 > y_0$. Then
\begin{equation}\label{1dim2}
x \in (-r)^n\tau_r^{-n}(x_0),\ y \in (-r)^n\tau_r^{-n}(y_0) \quad \mbox{implies that}\quad x > y.
\end{equation}
Since
\[
\mathcal{T}_r(z) = \mathop{\rm
Lim}_{n\rightarrow\infty} (-r)^n
\tau_r^{-n}(z)
\]
holds for each $z\in \mathbb{Z}$, \eqref{1dim2} yields that $\mathcal{T}_r(x_0)$ and $\mathcal{T}_r(y_0)$ have at most one point in common. Thus the result follows for the case $r>0$. The case  $r<0$ can be treated similarly.
\end{proof}

\section{SRS and canonical number systems}\label{cns}

The aim of this section is to relate SRS tiles to tiles
associated with expanding polynomials. We recall that an {\em expanding polynomial} is a polynomial each of whose roots is strictly larger than one in modulus.

\subsection{Expanding polynomials over $\mathbb{Z}$, SRS representations and canonical
number systems}

Let $A=a_d x^d+a_{d-1}x^{d-1}+\cdots+a_1 x+a_0 \in \mathbb{Z}[x]$, $a_0 \geq 2$, $a_d\neq 0$, and $\mathcal{Q}:=\mathbb{Z}[x]/A\mathbb{Z}[x]$ the factor ring, with $X \in \mathcal{Q}$ being the image of $x$ under the
canonical epimorphism. Furthermore, set
$\mathcal{N}=\{0,\ldots,a_0-1\}$.
We want to represent each element $P\in \mathcal{Q}$ formally as
\begin{equation}\label{cnsrepinfinity}
P = \sum_{n=0}^\infty b_n X^n \qquad(b_n \in \mathcal{N}).
\end{equation}
More precisely, we want to find a sequence $(b_n)_{n\in\mathbb{N}}$ as in the following definition.

\begin{definition}[$(A,\mathcal{N})$-representation]
Let $P \in \mathcal{Q}=\mathbb{Z}[x]/A\mathbb{Z}[x]$. A~sequence
$(b_n)_{n\in\mathbb{N}}$ with $b_n\in\mathcal{N}$ is called an
\emph{$(A,\mathcal{N})$-representation} of $P$ if $P -
\sum_{i=0}^{m-1} b_n X^n \in X^m\mathcal{Q}$ for all
$m\in\mathbb{N}$.
\end{definition}

In order to show that such an representation exists and is unique, we introduce a \emph{backward division mapping} $D_A:\mathcal{Q} \rightarrow \mathcal{Q}$.

\begin{lemma} \label{lem:bdm}
The backward division mapping $D_A:\mathcal{Q} \rightarrow \mathcal{Q}$ given for $P=\sum_{i=0}^\ell p_i X^i$, $p_i\in\mathbb{Z}$, by
\[
D_A(P) = \sum_{i=0}^{\ell-1} p_{i+1}X^i - \sum_{i=0}^{d-1} q a_{i+1}X^i,\quad
q=\left\lfloor{\frac{p_0}{a_0}}\right\rfloor,
\]
is well defined.
Every $P \in \mathcal{Q}$ has one and only one $(A,\mathcal{N})$-representation.
\end{lemma}

\begin{proof}
It is easy to see that $D_A(P)$ does not depend on the choice of the representation of~$P$, that $P=(p_0-q a_0)+X D_A(P)$, and that $b_0=p_0-q a_0$ is the unique element in $\mathcal{N}$ with $P-b_0 \in X\mathcal{Q}$ (for the case of monic polynomials $A$, this is detailed in \cite{Akiyama-Borbeli-Brunotte-Pethoe-Thuswaldner:05}; for the non-monic case, see \cite{Scheicher-Surer-Thuswaldner-vdWoestijne:08}).

The $(A,\mathcal{N})$-representation can be obtained by successively applying $D_A$, which yields
\[
P=\sum_{n=0}^{m-1} b_n X^n + X^m D_A^m(P)
\]
with $b_n=D_A^n(P)-X D_A^{n+1}(P)\in\mathcal{N}$.
Therefore, the $(A,\mathcal{N})$-representation of $P$ is unique.
\end{proof}

In order to relate $D_A$ to an SRS, we use an appropriate $\mathbb{Z}$-submodule of $\mathcal{Q}$, following \cite{Scheicher-Surer-Thuswaldner-vdWoestijne:08}.

\begin{definition}[Brunotte basis and Brunotte module]
The \emph{Brunotte basis modulo $A$} is defined by $\{W_0,\ldots,W_{d-1}\}$ with
\begin{equation} \label{eq:Brunottebasis}
W_0=a_d \quad\mbox{and}\quad W_k = X W_{k-1} + a_{d-k} \quad\mbox{for }1\le k\le d-1.
\end{equation}
The \emph{Brunotte module} $\Lambda_A$ is the $\mathbb{Z}$-submodule of $\mathcal{Q}$ generated by the Brunotte basis.
The representation mapping with respect to the Brunotte basis is denoted by
\[
\Psi_A: \,\Lambda_A \to \mathbb{Z}^d, \quad  P = \sum_{k=0}^{d-1} z_kW_k
\mapsto (z_0,\ldots,z_{d-1})^t.
\]
\end{definition}

Note that $\Lambda_A$ is isomorphic to $\mathcal{Q}$ if $A$ is monic. (Here, monic means that $a_d \in \{-1,1\}$.)
This is easily seen by checking that the coordinate matrix of $\{W_0,\ldots,W_{d-1}\}$ w.r.t. $\{1,X,\ldots,X^{d-1}\}$ is given by
\begin{equation}\label{eq:V}
V:=\left(\begin{array}{ccccc}
a_d & a_{d-1} & \cdots & \cdots &  a_1 \\
0 & \ddots & \ddots & & \vdots \\
\vdots & \ddots & \ddots & \ddots & \vdots \\
\vdots &  & \ddots & \ddots & a_{d-1} \\
0 & \cdots & \cdots & 0 & a_d
\end{array}\right).
\end{equation}

The restriction to Brunotte module $\Lambda_A$ allows us to relate the SRS transformation $\tau_\mathbf{r}$ to the backward division mapping $D_A$ in the following way.

\begin{proposition}\label{p:CNSconjugacy}
Let $A=a_d x^d+a_{d-1}x^{d-1}+\cdots+a_1x+a_0\in\mathbb{Z}[x]$, $a_0 \geq 2$, $a_d \neq 0$, $\mathbf{r}=\big(\frac{a_d}{a_0},\ldots,\frac{a_1}{a_0}\big)$.
Then we have
\begin{equation} \label{eq:CNSconj}
D_A^n\Psi_A^{-1}(\mathbf{z}) = \Psi_A^{-1}\tau_\mathbf{r}^n(\mathbf{z}) \quad\mbox{for all }\mathbf{z}\in\mathbb{Z}^d,\ n\in\mathbb{N}.
\end{equation}
In particular, the restriction of $D_A$ to the Brunotte module $\Lambda_A$ is conjugate to $\tau_\mathbf{r}$.
\end{proposition}

\begin{proof}
On $\Lambda_A$, the mapping $D_A$ can be written as
\begin{equation} \label{eq:TA}
D_A\bigg(\sum_{k=0}^{d-1} z_kW_k\bigg) = \sum_{k=0}^{d-2} z_{k+1}W_k - \left\lfloor \frac{z_0a_d+\cdots+z_{d-1}a_1}{a_0} \right\rfloor W_{d-1}
\end{equation}
(see \emph{e.g.}\ \cite[Section~3]{Akiyama-Borbeli-Brunotte-Pethoe-Thuswaldner:05}).
Therefore, we have $D_A\Psi_A^{-1}(\mathbf{z}) = \Psi_A^{-1}\tau_\mathbf{r}(\mathbf{z})$, which implies (\ref{eq:CNSconj}).
Since $\Psi_A$ is bijective, the restriction of $D_A$ to $\Lambda_A$ is conjugate to~$\tau_\mathbf{r}$.
\end{proof}

Note that in the monic case this gives a conjugacy between the backward division mapping on the full set $\mathcal{Q}$ and the SRS transformation.

With help of the conjugacy proved in
Proposition~\ref{p:CNSconjugacy}, we get a simple formula to gain
the $(A,\mathcal{N})$-representation (\ref{cnsrepinfinity}) for each
$P \in \Lambda_A$ using the associated
transformation~$\tau_\mathbf{r}$.

\begin{lemma}\label{lem:biCNS}
Let $A=a_d x^d+a_{d-1}x^{d-1}+\cdots+a_1x+a_0\in\mathbb{Z}[x]$, $a_0 \geq 2$, $a_d \neq 0$, $\mathbf{r}=\big(\frac{a_d}{a_0},\ldots,\frac{a_1}{a_0}\big)$.
The $(A,\mathcal{N})$-representation of $P \in \Lambda_A$ is given by
\[
b_n=\big\{\mathbf{r}\,\tau_\mathbf{r}^n\Psi_A(P)\big\}a_0 \quad\mbox{for all }n\in\mathbb{N}.
\]
\end{lemma}

\begin{proof}
For fixed $n$, let $D_A^n(P)=\sum_{k=0}^{d-1} z_k W_k$ and recall that $\Psi_AD_A^n(P)=(z_0,\ldots,z_{d-1})$.
By (\ref{eq:Brunottebasis}), (\ref{eq:TA}) and the fact that $X W_{d-1}+a_0=0$, we obtain that
\[
X D_A\bigg(\sum_{k=0}^{d-1} z_kW_k\bigg) = \sum_{k=0}^{d-2} z_{k+1}(W_{k+1}-a_{d-k-1}) + \left\lfloor \frac{z_0a_d+\cdots+z_{d-1}a_1}{a_0} \right\rfloor a_0,
\]
therefore
\begin{align*}
b_n & = D_A^n(P)-X D_A^{n+1}(P) \\
& = z_0W_0 + \sum_{k=0}^{d-2}z_{k+1}a_{d-k-1} - \left\lfloor \frac{z_0a_d+\cdots+z_{d-1}a_1}{a_0} \right\rfloor a_0 = \left\{ \frac{z_0a_d+\cdots+z_{d-1}a_1}{a_0} \right\} a_0 \\
& = \big\{\mathbf{r}\,\Psi_A D_A^n(P)\big\} a_0 = \big\{\mathbf{r}\,\tau_\mathbf{r}^n\Psi_A(P)\big\} a_0. \hfill\qedhere
\end{align*}
\end{proof}

If the $(A,\mathcal{N})$-representation $(b_n)_{n\in\mathbb{N}}$
has only finitely many non-zero elements, then $P$ can be written as a finite sum of the shape
\[
P=\sum_{n=0}^{m-1} b_n X^n \qquad (b_n \in \mathcal{N}).
\]
We recover the following well-known notion of canonical number systems.

\begin{definition}[Canonical number system]
Let $A=a_d x^d+\cdots+a_1x+a_0 \in \mathbb{Z}[x]$, $a_0\ge2$, $a_d \neq 0$, $\mathcal{Q}=\mathbb{Z}[x]/A\mathbb{Z}[x]$, and
$\mathcal{N}=\{0,\ldots,a_0-1\}$.
If for each $P \in \mathcal{Q}$, the $(A,\mathcal{N})$-representation $(b_n)_{n\in\mathbb{N}}$, $b_n\in\mathcal{N}$, has only finitely many non-zero elements, then we call $(A,\mathcal{N})$ a \emph{canonical number system} (CNS, for short).
\end{definition}

It is shown in \cite{Scheicher-Surer-Thuswaldner-vdWoestijne:08}
that it is sufficient to check the finiteness of the
$(A,\mathcal{N})$-representations for all $P\in\Lambda_A$ in order
to check whether $(A,\mathcal{N})$ is a CNS. Thus the conjugacy
between $D_A$ and $\tau_\mathbf{r}$ is sufficient to reformulate the
CNS property in terms of SRS with finiteness property.

\begin{proposition} \label{prop:CNS}
Let $A=a_d x^d+a_{d-1}x^{d-1}+\cdots+a_1 x+a_0\in\mathbb{Z}[x]$, $a_0 \geq 2$, $a_d \neq 0$, and $\mathcal{N}=\{0,\ldots,a_0-1\}$.
Then the following assertions hold.
\begin{itemize}
\item
The polynomial $A$ is expanding if and only if $\mathbf{r}=\big(\frac{a_d}{a_0},\frac{a_{d-1}}{a_0},\ldots,\frac{a_1}{a_0}\big)$ is contained in $\mathrm{int}(\mathcal{D}_d)$.
\item
The pair $(A,\mathcal{N})$ is a CNS if and only if $\mathbf{r}=\big(\frac{a_d}{a_0},\frac{a_{d-1}}{a_0},\ldots,\frac{a_1}{a_0}\big) \in \mathcal{D}_d^{(0)}$.
\end{itemize}
\end{proposition}

\begin{proof}
The first assertion follows from Lemma~\ref{lem:int}.
The second assertion follows from the conjugacy given in  Proposition~\ref{p:CNSconjugacy} together with the fact that it is sufficient to check the finiteness of the $(A,\mathcal{N})$-representations for all $P\in\Lambda_A$ (see also \cite{Scheicher-Surer-Thuswaldner-vdWoestijne:08}).
\end{proof}

\subsection{Tiles associated with an expanding polynomial}
We define two classes of tiles for expanding polynomials.
The first definition goes back to K\'atai and K\H{o}rnyei~\cite{Katai-Koernyei:92} (see also~\cite{Scheicher-Thuswaldner:01}) for monic polynomials which give rise to a CNS.
We first extend this definition to arbitrary expanding polynomials.

\begin{definition}[Self-affine tile associated with an expanding polynomial]\label{CNStiledef}
Let $A=a_d x^d+a_{d-1}x^{d-1}+\cdots+a_1 x+a_0\in \mathbb{Z}[x]$ be an expanding polynomial with $a_0\ge 2$, $a_d \neq 0$, and  $\mathcal{N}=\{0,\ldots,a_0-1\}$.
The \emph{self-affine tile associated with $A$} is defined as
\begin{equation}\label{CNStile}
\mathcal{F}_A :=\Big\{ \mathbf{t} \in \mathbb{R}^d \;\Big|\;
\mathbf{t}=\sum_{i=1}^\infty B^{-i}(c_i,0,\ldots,0)^t,\ c_i \in
\mathcal{N}\Big\}
\end{equation}
where $B:=V M_\mathbf{r}^{-1}V^{-1}$ with $\mathbf{r}=\big(\frac{a_d}{a_0},\frac{a_{d-1}}{a_0},\ldots,\frac{a_1}{a_0}\big)$, and $V$ is given by (\ref{eq:V}).
\end{definition}

\begin{remark}
Easy calculations show that the matrix $B$ in Definition~\ref{CNStiledef} is given by
\[
B =  \begin{pmatrix}
0 & \cdots & \cdots & 0 & -\frac{a_0}{a_d} \\
1 & \ddots & & \vdots  & -\frac{a_1}{a_d} \\
0 & \ddots & \ddots & \vdots & \vdots \\
\vdots & \ddots & \ddots & 0 & -\frac{a_{d-2}}{a_d} \\
\vphantom{\vdots}0 & \cdots & 0 & 1 & -\frac{a_{d-1}}{a_d}
\end{pmatrix}.
\]
We use it instead of $M_\mathbf{r}^{-1}$ in order to be consistent with the existing literature on CNS tiles.
\end{remark}

Since $A$ is expanding, it is easy to see that the series in
(\ref{CNStile}) always converges.
The tile $\mathcal{F}_A$ has the following properties.

\begin{itemize}
\item
The tile $\mathcal{F}_A$ is compact.
\item
The tile $\mathcal{F}_A$ is a self-affine set (hence the terminology) as it obeys the set equation
$B \mathcal{F}_A = \bigcup_{c \in \mathcal{N}} \big(\mathcal{F}_A +
(c,0\ldots,0)^t\big)$.
Indeed, it is the unique non-empty compact set satisfying this equation (\emph{cf.\ e.g.}\ \cite{Barnsley:88,Falconer:90,Hutchinson:81}).
Self-affine tiles have been studied extensively in a very general
context in the literature.
We refer the reader to the surveys by Vince~\cite{Vince:00a} and Wang~\cite{Wang:99}.
\item
If $A$ is monic, then $\{\mathbf{z}+\mathcal{F}_A \mid \mathbf{z}\in\mathbb{Z}^d\}$ forms a tiling of $\mathbb{R}^d$ if $A$ is an irreducible polynomial (this is an immediate consequence of \cite[Corollary~6.2]{Lagarias-Wang:97}). Moreover, there exist algorithms in order to decide whether $\{\mathbf{z}+\mathcal{F}_A \mid \mathbf{z}\in\mathbb{Z}^d\}$ forms a tiling for any given polynomial $A$ (see \emph{e.g.}\ \cite{Vince:00a}).
\end{itemize}

For non-monic polynomials $A$, the above definition turns out to be
not well-suited. The tiles have strong overlaps. The submodule
$\Lambda_A$ is a good tool to define a new class of tiles which
forms a (multiple) tiling also for non-monic polynomials.

\begin{definition}[Brunotte tile associated with an expanding polynomial]\label{Brunottetiledef}
Let $A =a_d x^d+a_{d-1}x^{d-1}+\cdots+a_1 x+a_0\in \mathbb{Z}[x]$ be an expanding polynomial with $a_0\ge 2$, $a_d \neq 0$, and $\mathcal{N}=\{0,\ldots,a_0-1\}$.
For each $P\in\Lambda_A$, the \emph{Brunotte tile} associated with $A$ is defined as
\[
\mathcal{G}_A(P) := \mathop{\rm Li}_{n\to\infty} B^{-n}V\Psi_A\big(D_A^{-n}(P) \cap \Lambda_A\big),
\]
where $\mathop{\rm Li}$ denotes the lower Hausdorff limit.
The set $\mathcal{G}_A(0)$ is called the \emph{central Brunotte tile} associated with the expanding polynomial~$A$.
\end{definition}

\begin{lemma} \label{lem:FLim}
We have
\[
\mathcal{F}_A = \mathop{\rm Lim}_{n\to\infty} B^{-n}V\Psi_A\big(D_A^{-n}(0)\big)
\]
and thus $\mathcal{G}_A(P) \subseteq V\Psi_A(P) + \mathcal{F}_A$.
\end{lemma}

\begin{proof}
By the proof of Lemma~\ref{lem:bdm}, we have $P \in D_A^{-n}(0)$ if and only if $P = \sum_{i=1}^n c_i X^{n-i}$ for some $c_i \in \mathcal{N}$.
Since $V \Psi_A(P) = (p_0,\ldots,p_{d-1})^t$ if $P = \sum_{i=0}^{d-1} p_i X^i$, we obtain recursively that $V \Psi_A(X^i) = B^i (1,0,\ldots,0)^t$ for all $i \ge 0$, hence $B^{-n}V\Psi_A(\sum_{i=1}^n c_i X^i) = \sum_{i=1}^n B^{-i}(c_i,0,\ldots,0)^t$, which yields the lemma.
\end{proof}

If $A$ is monic, then we will show that the inclusion in Lemma~\ref{lem:FLim} becomes the equality $\mathcal{G}_A(P)=V\Psi_A(P)+\mathcal{F}_A$ for every $P\in\mathcal{Q}$.

\subsection{From tiles associated with expanding polynomials to SRS tiles}

Next we show that the Brunotte tiles are obtained from the SRS tiles by a linear transformation.

\begin{theorem}\label{thm:FT}
Let $A=a_d x^d+a_{d-1}x^{d-1}+\cdots+a_1 x+a_0 \in \mathbb{Z}[x]$, $a_0 \geq 2$, $a_d \neq 0$, be an expanding polynomial, and $\mathbf{r}=\big(\frac{a_d}{a_0},\frac{a_{d-1}}{a_0},\ldots,\frac{a_1}{a_0}\big)$.
Then
\[
\mathcal{G}_A\big(\Psi_A^{-1}(\mathbf{z})\big) = V \mathcal{T}_\mathbf{r}(\mathbf{z}) \quad\mbox{for all }\mathbf{z}\in\mathbb{Z}^d.
\]
\end{theorem}

\begin{proof}
Set $P=\Psi_A^{-1}(\mathbf{z})$.
We have
\begin{equation} \label{eq:Brunotteapprox}
B^{-n}V\Psi_A(D_A^{-n}(P)\cap\Lambda_A) = B^{-n}V\tau_\mathbf{r}^{-n}(\mathbf{z}) = V M_\mathbf{r}^n\tau_\mathbf{r}^{-n}(\mathbf{z}).
\end{equation}
Taking $\mathop{\rm Li}$ for $n\to\infty$ proves the result.
\end{proof}

The following corollary shows that the lower Hausdorff limit in the definition of the Brunotte tiles is indeed a Hausdorff limit, as for the SRS tiles.
Furthermore, we show how (\ref{CNStile}) has to be adapted to give the Brunotte tiles.

\begin{corollary}
The Brunotte tiles can be written as
\[
\mathcal{G}_A(P) = \mathop{\rm Lim}_{n\to\infty} B^{-n}V\Psi_A\big(D_A^{-n}(P) \cap \Lambda_A\big),
\]
where $\mathop{\rm Lim}$ denotes the Hausdorff limit.
Moreover, if $P=\sum_{k=0}^{d-1}p_k X^k$, then
\begin{align*}
\mathcal{G}_A(P) = (p_0,\ldots,p_{d-1})^t + \Big\{ \mathbf{t} \in \mathbb{R}^d \;\Big|\; &
\mathbf{t}=\sum_{i=1}^\infty B^{-i}(c_i,0,\ldots,0)^t,\ c_i \in
\mathcal{N}, \\
& \sum_{i=1}^n c_i X^{n-i} + X^n P \in \Lambda_A\ \mbox{for all}\ n\in\mathbb{N}\Big\}.
\end{align*}
\end{corollary}

\begin{proof}
The first assertion follows from (\ref{eq:Brunotteapprox}) and Theorem~\ref{seteq}.
The second assertion follows from Proposition~\ref{prop:rep} and the fact that
\[
D_A^{-n}(P) \cap \Lambda_A = \bigg\{\sum_{i=1}^n c_i X^{n-i} + X^n P \in \Lambda_A \;\bigg|\; c_i\in\mathcal{N}\bigg\}. \hfill\qedhere
\]
\end{proof}

For the monic case, we derive the following corollary.

\begin{corollary}\label{2904083}
Suppose that the polynomial $A$ in Theorem~\ref{thm:FT} is monic.
Then
\[
\mathcal{G}_A(P) = V\Psi_A(P)+\mathcal{F}_A \quad\mbox{for all }P\in\mathcal{Q}.
\]
\end{corollary}

Now we can state the main result of the present section.
It establishes tiling properties for Brunotte tiles which are even valid in the non-monic case.

\begin{theorem}
Let $A=a_d x^d+a_{d-1}x^{d-1}+\cdots+a_1 x+a_0 \in \mathbb{Z}[x]$, $a_0 \geq 2$, $a_d \neq 0$, be an expanding polynomial, and $\mathcal{N}=\{0,\ldots,a_0-1\}$.
Then the following assertions hold.
\begin{itemize}
\item
The collection $\{\mathcal{G}_A(P) \mid P\in\Lambda_A\}$ forms a weak $m$-tiling of $\mathbb{R}^d$ for some $m\ge1$.
\item
If $(A,\mathcal{N})$ is a CNS, then $\{\mathcal{G}_A(P) \mid P\in\Lambda_A\}$ forms a weak tiling of $\mathbb{R}^d$.
\end{itemize}
\end{theorem}

\begin{proof}
By Theorem~\ref{thm:FT}, it is equivalent to consider the collection $\{\mathcal{T}_\mathbf{r}(\mathbf{z}) \mid \mathbf{z}\in\mathbb{Z}^d\}$ for $\mathbf{r}=\big(\frac{a_d}{a_0},\frac{a_{d-1}}{a_0},\ldots,\frac{a_1}{a_0}\big)$.
In view of Proposition~\ref{prop:CNS}, the first assertion follows from Theorem~\ref{theo:mdense}, while the second assertion is a consequence of Corollary~\ref{cor:weaktiling}.
\end{proof}

\begin{example}
Consider the monic CNS polynomial $A=x^2-x+2$. The associated SRS parameter is $\mathbf{r}=(\frac{1}{2},-\frac{1}{2})$.
The central SRS tile $\mathcal{T}_\mathbf{r}(\mathbf{0})$ as well as its neighbors are shown in Figure~\ref{tm12} on the left hand side.
To obtain $\mathcal{G}_A\big(\Psi^{-1}(\mathbf{z})\big) = V\mathbf{z}+\mathcal{F}_A$, we have to multiply the SRS tiles by the matrix $V=\left(\begin{array}{cc} 1& -1 \\ 0 & 1\end{array}\right)$.
The Brunotte tiles are shown on the right hand side of Figure~\ref{tm12}.
\begin{figure}[h]
\centering
\includegraphics[width=.35\textwidth]{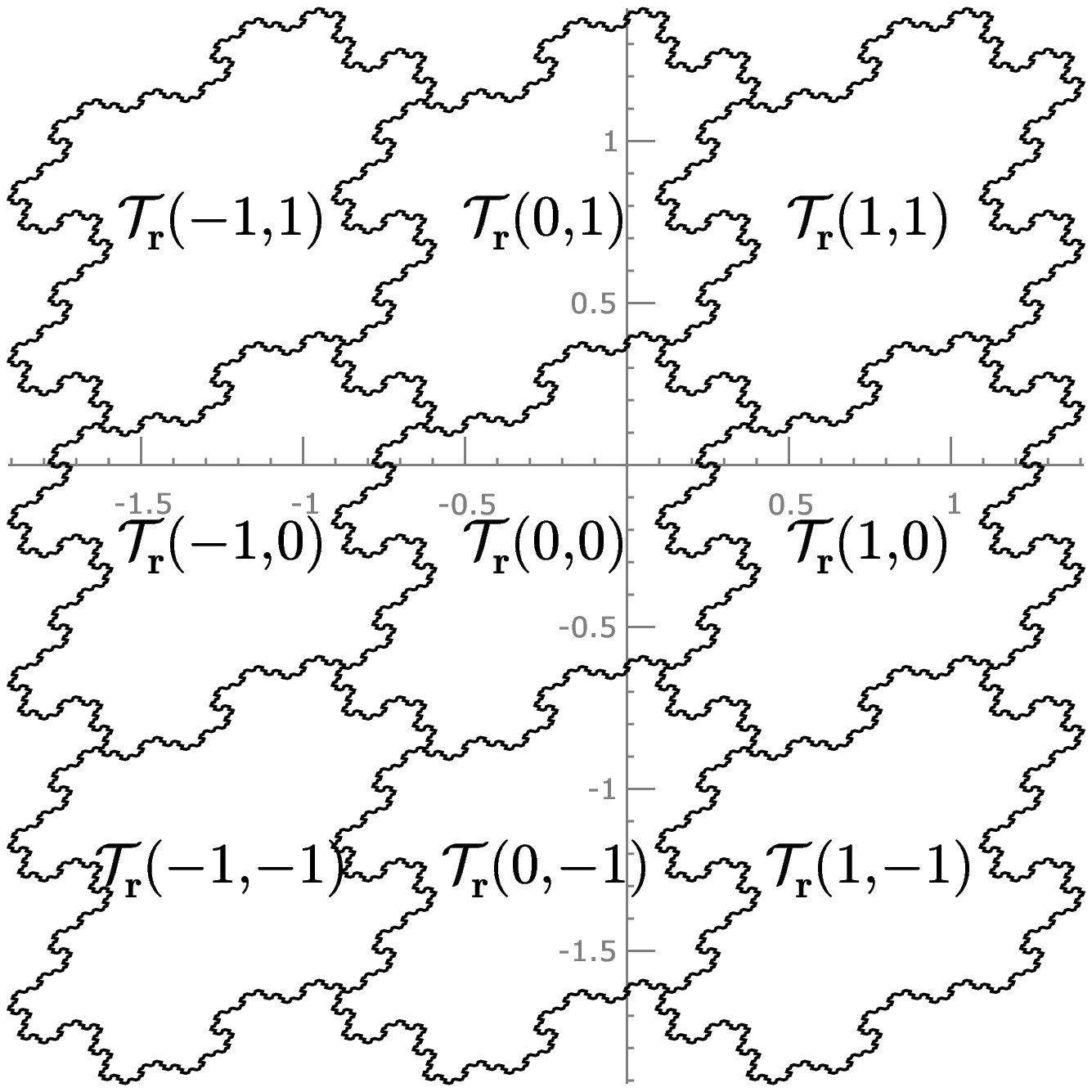}
\quad
\includegraphics[width=.5925\textwidth]{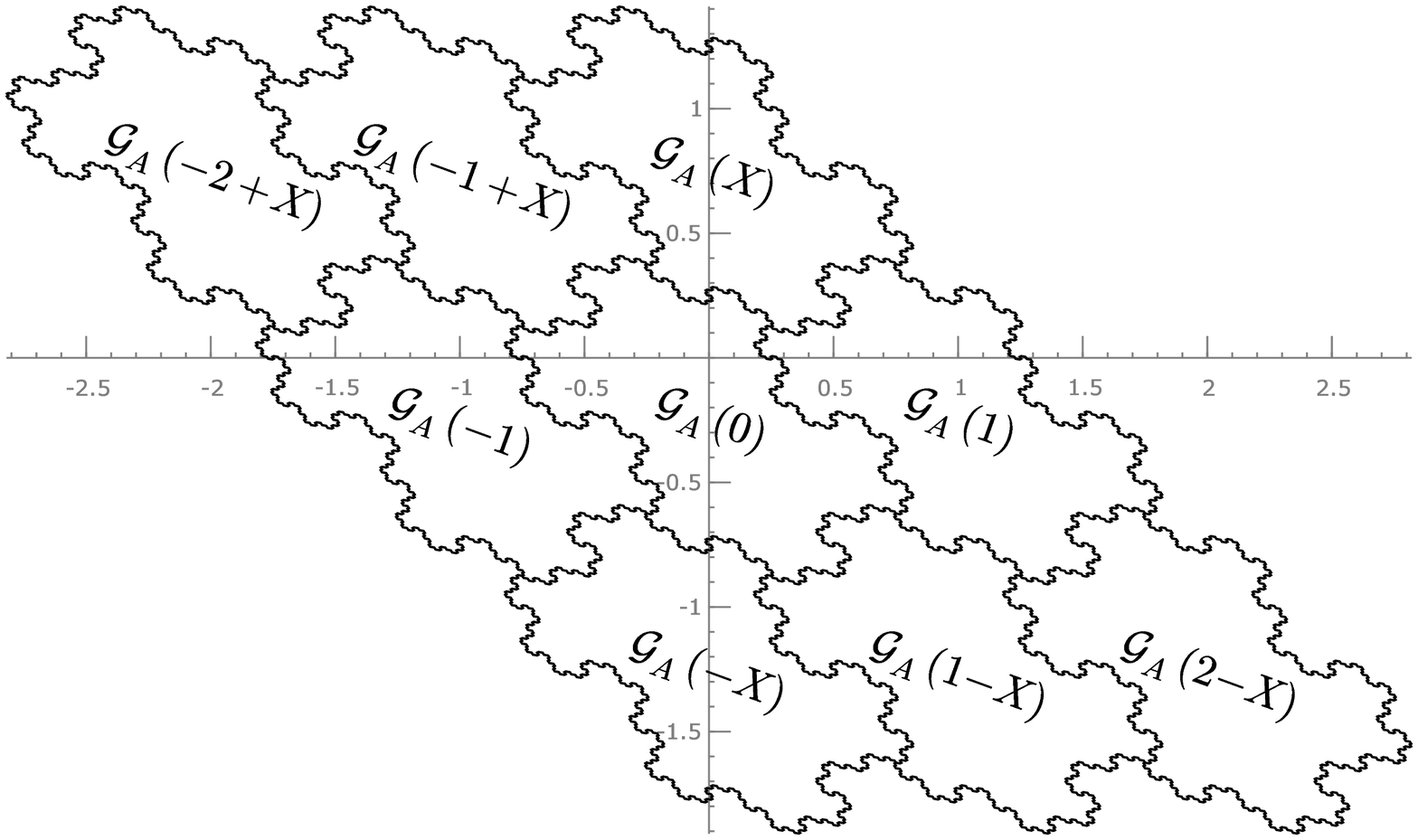}
\caption{The SRS tiles $\mathcal{T}_\mathbf{r}(\mathbf{z})$ for $\mathbf{r}=\left(\frac{1}{2},-\frac{1}{2}\right)$, $\|\mathbf{z}\|_\infty\le 1$, (left) and the Brunotte tiles $\mathcal{G}_A\big(\Psi^{-1}(\mathbf{z})\big) = V\mathbf{z}+\mathcal{F}_A$ associated with $A=x^2-x+2$ (right).} \label{tm12}
\end{figure}
\end{example}

An example of SRS tiles associated with the parameter
$\mathbf{r}=(3/4,1)$ is discussed in Example~\ref{ex:bt}. The tiles corresponding to this parameter are depicted in Figure~\ref{map}. It is related to a non-monic CNS. Also the parameter $\mathbf{r}=(\frac{9}{10},-\frac{11}{20})$, which yields tiles that are single points, corresponds to a non-monic CNS (see Example~\ref{ex:pp}).

\subsection{Rational base number systems}\label{rbns}

Akiyama \emph{et al.}~\cite{Akiyama-Frougny-Sakarovitch:07}
considered expansions of integers in rational bases $p/q$, with
coprime integers $p>q\ge1$, of the form
\[
N = \frac{1}{q} \sum_{n=0}^\infty b_n \Big(\frac{p}{q}\Big)^n \qquad(b_n\in\mathcal{N}=\{0,\ldots,p-1\}).
\]
In our setting, the sequence $(b_n)_{n\in\mathbb{N}}$ is the
$(A,\mathcal{N})$-representation of $q N$, where $A=-q x+p$. The
Brunotte basis modulo $A$ is given by $\{-q\}$. By
Lemma~\ref{lem:biCNS}, the corresponding SRS $\tau_{-q/p}(N) =
-\big\lfloor -N\frac{q}{p} \big\rfloor$ yields
$b_n=\big\{-\frac{q}{p} \tau_{-q/p}^n(-N)\big\}p$. (Here we write
one-dimensional vectors as scalars.) In view of Theorem~\ref{d1}, we
obtain that the collection of tiles associated with these number
systems forms a tiling for each choice $p/q$ with $p>q\ge1$.

\medskip

For the case $p/q=3/2$, we show that the collection
$\{\mathcal{T}_{-2/3}(N) \mid N\in\mathbb{Z}\}$ consists of
(possibly degenerate) intervals with infinitely many different
lengths.

\begin{lemma} \label{lem:513}
Let $I$ be a finite set of consecutive integers.
Then $\tau_{-2/3}^{-n}(I)$ is a finite set of consecutive integers for all $n\in\mathbb{N}$, and we have
\[
(\#I-1)\Big(\frac32\Big)^n + 1 \le \#\tau_{-2/3}^{-n}(I) \le (\#I+1)\Big(\frac32\Big)^n - 1.
\]
\end{lemma}

\begin{proof}
By the definition of $\tau_{-2/3}$, the preimage of every finite set
of consecutive integers is a finite set of consecutive integers (see
also the proof of Theorem~\ref{d1}). We have
$\#\tau_{-2/3}^{-1}(N)=2$ if and only if $N$ is even and
$\#\tau_{-2/3}^{-1}(N)=1$ if $N$ is odd. Therefore the inequalities
hold for $n=1$, and by induction for all $n\in\mathbb{N}$.
\end{proof}

\begin{lemma} \label{lem:514}
For every $k\ge1$, there exists some $N_k\in\mathbb{Z}$ such that $\#\tau_{-2/3}^{-k}(N_k)=2$.
\end{lemma}

\begin{proof}
It follows from
\cite[Proposition~10]{Akiyama-Frougny-Sakarovitch:07} that there
exists some $L_k\in\mathbb{Z}$ such that the
$(-2x+3,\{0,1,2\})$-representation of $2L_k$ satisfies $b_0=0$,
$b_n=1$ for $1\le n<k$. The SRS representation of $-L_k$ with
respect to $-2/3$ is thus given by
$(b_0/3,b_1/3,\ldots)=(0,1/3,\ldots,1/3,b_k/3,b_{k+1}/3,\ldots)$.
Let $N_k=\tau_{-2/3}^k(-L_k)$. For $1\le n<k$, the set
$\tau_{-2/3}^{-n}(N_k)$ consists only of the number with SRS
representation $(1/3,\ldots,1/3,b_k/3,b_{k+1}/3,\ldots)$ because
$\big(1/3+(-2/3)\mathbb{Z}\big)\cap[0,1)=\{1/3\}$. Thus
$\tau_{-2/3}^{-k}(N_k)$ consists exactly of the two numbers $-L_k$
and $-L_k-1$, which have the SRS representations
$(0,1/3,\ldots,1/3,b_k/3,b_{k+1}/3,\ldots)$ and
$(2/3,1/3,\ldots,1/3,b_k/3,b_{k+1}/3,\ldots)$, respectively.
\end{proof}

\begin{proposition}
For every $k\ge1$, there exists some $N_k\in\mathbb{Z}$ such that
\[
\Big(\frac32\Big)^{-k} \le \lambda_1(\mathcal{T}_{-2/3}(N_k)) \le 3\Big(\frac32\Big)^{-k},
\]
where $\lambda_1$ denotes the one-dimensional Lebesgue measure.
\end{proposition}

\begin{proof}
We have
\[
\lambda_1(\mathcal{T}_{-2/3}(N_k)) = \lim_{n\to\infty} \Big(\frac23\Big)^n \#\tau_{-2/3}^{-n}(N_k) = \lim_{n\to\infty} \Big(\frac23\Big)^n \#\tau_{-2/3}^{k-n}(\tau_{-2/3}^{-k}(N_k))
\]
for all $N_k\in\mathbb{Z}$.
For the $N_k$ given by Lemma~\ref{lem:514}, we have $\#\tau_{-2/3}^{-k}(N_k)=2$.
Using Lemma~\ref{lem:513}, the inequalities are proved.
\end{proof}

Summing up we get the following result.

\begin{corollary} \label{cor:515}
The tiling $\{\mathcal{T}_{-2/3}(N) \mid N\in\mathbb{Z}\}$ consists
of (possibly degenerate) intervals with infinitely many different
lengths.
\end{corollary}

\section{SRS and beta-expansion}\label{beta}
We prove now that there is a tight relation with beta-tiles when
$\mathbf{r}$ is related to a unit Pisot number, and that SRS tiles are new
objects in the non-unit case (similarly as shown in
Section~\ref{cns} in the non-monic CNS case).

\subsection{Beta-expansions and SRS representations}\label{Tbeta}
We start with beta-expansions, which were first studied by
R{\'e}nyi~\cite{Renyi:57} and Parry~\cite{Parry:60}. For a real
number $\beta>1$, the \emph{$\beta$-transformation}
$T_\beta:[0,1)\to[0,1)$ is defined by $T_\beta(x)=\{\beta x\}=\beta
x-\lfloor\beta x\rfloor$. The \emph{$\beta$-expansion} of
$x\in[0,1)$ is
\[
x=\sum_{n=1}^\infty b_n \beta^{-n} \quad\mbox{with}\quad b_n= \big\lfloor\beta T_\beta^{n-1}(x)\big\rfloor\quad \mbox{for all}\ n\ge1.
\]
The following relation between $T_\beta$ and $\tau_\mathbf{r}$ was
shown in \cite{Akiyama-Borbeli-Brunotte-Pethoe-Thuswaldner:05}, see
also~\cite{Hollander:96}.

\begin{proposition}\label{prop:betanumformula}
Let $\beta>1$ be an algebraic integer with minimal polynomial
\begin{equation} \label{er}
x^{d+1}+a_d x^d+\cdots+a_1x+a_0 =
(x-\beta)(x^d+r_{d-1}x^{d-1}+\cdots+r_1x+r_0)
\end{equation}
and $\mathbf{r}=(r_0,\ldots,r_{d-1})$.
Then we have
\begin{equation} \label{eq:conj}
\{\mathbf{r}\tau_\mathbf{r}^n(\mathbf{z})\}=T_\beta^n(\{\mathbf{r}\mathbf{z}\}) \quad \mbox{for all } \mathbf{z}\in\mathbb{Z}^d,\,n\in\mathbb{N}.
\end{equation}
In particular, the restriction of $T_\beta$ to $\mathbb{Z}[\beta]\cap[0,1)$ is conjugate to $\tau_\mathbf{r}$.
\end{proposition}

\begin{proof}
Let $\mathbf{z}=(z_0,\ldots,z_{d-1})^t\in\mathbb{Z}^d$.
If we set $z_d=-\lfloor\mathbf{r}\mathbf{z}\rfloor$, then we have
\begin{equation} \label{eq:rz}
\{\mathbf{r}\mathbf{z}\}=(r_0,\ldots,r_{d-1},1)(z_0,\ldots,z_{d-1},z_d)^t.
\end{equation}
Furthermore, $(r_0,\ldots,r_{d-1},1)$ is a left eigenvector of the
companion matrix $M_{(a_0,\ldots,a_d)}$, in particular
\begin{equation} \label{eq:rMa}
(r_0,\ldots,r_{d-1},1)M_{(a_0,\ldots,a_d)}=\beta(r_0,\ldots,r_{d-1},1).
\end{equation}
Using (\ref{eq:rz}), (\ref{eq:rMa}) and the fact that $M_{(a_0,\ldots,a_d)}(z_0,\ldots,z_d)^t = (z_1,\ldots,z_d,m)^t$ with $m\in\mathbb{Z}$, we gain
\[
\{\mathbf{r}\tau_\mathbf{r}(\mathbf{z})\} = \{\mathbf{r}(z_1,\ldots,z_d)^t\} =  \{(r_0,\ldots,r_{d-1},1)M_{(a_0,\ldots,a_d)}(z_0,\ldots,z_d)^t\} = \{\beta\{\mathbf{r}\mathbf{z}\}\} = T_\beta(\{\mathbf{r}\mathbf{z}\}),
\]
Inductively, we obtain (\ref{eq:conj}).
Since the polynomial in (\ref{er}) is irreducible,  $\{r_0,\ldots,r_{d-1},1\}$ is a basis of $\mathbb{Z}[\beta]$.
Therefore the map $f:\mathbb{Z}^d\to\mathbb{Z}[\beta]\cap[0,1)$, $\mathbf{z}\mapsto\{\mathbf{r}\mathbf{z}\}$ is bijective, and we have $f\circ\tau_\mathbf{r}=T_\beta\circ f$.
\end{proof}

\begin{corollary} \label{cor:betanumformula}
Let $\beta$ and $\mathbf{r}$ be defined as in
Proposition~\ref{prop:betanumformula}, $(v_1,v_2,v_3,\ldots)$ be the
SRS representation of $\mathbf{z} \in \mathbb{Z}^d$ and
$\{\mathbf{r}\mathbf{z}\}=\sum_{n=1}^\infty b_n \beta^{-n}$  be the
$\beta$-expansion of $v_1=\{\mathbf{r}\mathbf{z}\}$. Then we have
\[
v_n=T_\beta^{n-1}(\{\mathbf{r}\mathbf{z}\}) \quad \mbox{and} \quad b_n=\beta v_n-v_{n+1} \quad \mbox{for all }n\ge1.
\]
\end{corollary}

\begin{proof}
By Definition~\ref{srsex} and (\ref{eq:conj}), we have $v_n=\{\mathbf{r}\tau_\mathbf{r}^{n-1}(\mathbf{z})\}=T_\beta^{n-1}(\{\mathbf{r}\mathbf{z}\})$, which yields the first assertion.
Using this equation and the definition of the $\beta$-expansion, we obtain
\[
b_n = \big\lfloor \beta T_\beta^{n-1}(\{\mathbf{r}\mathbf{z}\}) \big\rfloor = \beta T_\beta^{n-1}(\{\mathbf{r}\mathbf{z}\}) - \big\{\beta T_\beta^{n-1}(\{\mathbf{r}\mathbf{z}\})\big\} = \beta v_n - T_\beta^n(\{\mathbf{r}\mathbf{z}\}) = \beta v_n-v_{n+1}. \hfill\qedhere
\]
\end{proof}

The eigenvalues of $M_\mathbf{r}$ are exactly the Galois conjugates of $\beta$.
It follows by Lemma~\ref{lem:int} that $\mathbf{r} \in \mathrm{int}(\mathcal{D}_d)$ if $\beta$ is a Pisot number.

\begin{definition}[Finiteness property (F)]\label{def:F}
A number $\beta>1$ is said to have the \emph{finiteness property (F)} if the $\beta$-expansion of every $x \in \mathbb{Z}[\beta^{-1}] \cap [0,1)$ is finite, \emph{i.e.}, $T_\beta^n(x)=0$ for some $n\in\mathbb{N}$.
\end{definition}

Frougny and Solomyak~\cite{Frougny-Solomyak:92} proved that (F)
implies that $\beta$ is a Pisot number. By
\cite[Lemma~3]{Akiyama:02}, it is sufficient to consider $x \in
\mathbb{Z}[\beta] \cap [0,1)$ in Definition~\ref{def:F} (note that
in the present section $\mathbb{Z}[\beta] \cap [0,1)$ plays the same
role as $\Lambda_A$ plays in Section~\ref{cns}). Therefore
Proposition~\ref{prop:betanumformula} implies the following result
(see also
\cite[Theorem~2.1]{Akiyama-Borbeli-Brunotte-Pethoe-Thuswaldner:05}).

\begin{proposition}\label{prop:beta}
Let $\beta>1$ be an algebraic integer and let $\mathbf{r}$ be as in
Proposition~\ref{prop:betanumformula}. Then the following assertions
hold.
\begin{itemize}
\item
$\mathbf{r}\in\mathrm{int}(\mathcal{D}_d)$ if and only if $\beta$ is
a Pisot number.

\item
$\mathbf{r} \in \mathcal{D}_d^{(0)}$  if and only if $\beta$
satisfies (F).
\end{itemize}
\end{proposition}

\subsection{Beta-tiles}
For any Pisot number $\beta$ of degree $d+1$, Akiyama~\cite{Akiyama:02} defined a family of tiles covering $\mathbb{R}^d$ which is conjectured to be always a tiling if $\beta$ is a unit Pisot number, \emph{i.e.}, if $|a_0|=1$ in (\ref{er}).
If $\beta$ is not a unit, then this family cannot be a tiling. In analogy with the previous section, we modify the definition of the tiles so that we obtain tilings also for non-unit Pisot numbers.

Let $\beta_1,\ldots,\beta_d$ be the $d=r+2s$ Galois conjugates of $\beta$, such that $\beta_1,\ldots,\beta_r\in\mathbb{R}$ and $\beta_{r+1},\ldots,\beta_{r+2s} \in \mathbb{C}\setminus\mathbb{R}$ with $\beta_{r+1}=\overline{\beta_{r+s+1}},\,\ldots,\,\beta_{r+s}=\overline{\beta_{r+2s}}$.
For $x \in \mathbb{Q}(\beta)$, $1\le j\le d$, denote by $x^{(j)} \in \mathbb{Q}(\beta_j)$ the corresponding conjugate of $x$, and let
\[
\Phi_\beta:\, \mathbb{Q}(\beta) \to \mathbb{R}^d,\ x  \mapsto
\big(x^{(1)},\ldots,x^{(r)},\Re\big(x^{(r+1)}\big),\Im\big(x^{(r+1)}\big),\ldots,\Re\big(x^{(r+s)}\big),\Im\big(x^{(r+s)}\big)\big)^t.
\]

\begin{definition}[Beta-tile; see \cite{Akiyama:02}] \label{def:fracbeta}
Let $\beta$ be a Pisot number.
For $x \in \mathbb{Z}[\beta^{-1}] \cap [0,1)$, the set
\[
\mathcal{R}_\beta(x) := \mathop{\rm
Lim}_{n\to\infty}\Phi_\beta\big(\beta^n T_\beta^{-n}(x)\big),
\]
where the limit is taken with respect to the Hausdorff metric, is called a \emph{$\beta$-tile}.
The tile $\mathcal{R}_\beta(0)$ is called \emph{central $\beta$-tile}.
\end{definition}

Note that the limit in Definition~\ref{def:fracbeta} exists since
$\Phi_\beta\big(\beta^n T_\beta^{-n}(x)\big) \subseteq
\Phi_\beta\big(\beta^{n+1}T_\beta^{-n-1}(x)\big)$. Indeed, if
$y\in[0,1)$, then $T_\beta(y\beta^{-1})=y$, which implies $y\in\beta
T^{-1}(y)$.

For unit Pisot numbers, beta-tiles have been studied extensively (mostly under the name ``central tile'') (see \emph{e.g.}\ \cite{Akiyama:98,Akiyama:02,Thurston:89}) and are strongly related to Rauzy fractals associated with unit Pisot substitutions.
For a recent survey on these relations and on properties of Rauzy fractals and beta-tiles we refer to~\cite{Berthe-Siegel:05}.

Beta-tiles are known to satisfy a graph-directed IFS equation to be
compared with the set equation of Theorem~\ref{seteq} (see
\emph{e.g.}\ the survey \cite{Berthe-Siegel:05}). Nonetheless the
pieces obtained in the decomposition of the central tile might not
be measurably disjoint. When $\beta$ is a unit Pisot number, it is
proved that these pieces are disjoint, but when $\beta$ is not a
unit, the pieces do overlap (\emph{cf.}~\cite{Siegel:02}). To make
the pieces disjoint in the non-unit case, one can use two different
strategies: enlarging the representation space or making the tiles
smaller. The first strategy can be found in the literature (see
\cite{ABBS:08,Siegel:02}). It consists in adding $p$-adic fields for
the prime divisors $p$ of the norm of $\beta$ to the representation
space. In the present work, we want to carry out the second
strategy, leading to tiles that form a tiling of $\mathbb{R}^d$. To
this matter, in the $p$-adic extension obtained with the first
strategy, we have to choose among all points with the same Euclidean
part, a single specific point having $p$-adic components with
certain properties. Using Proposition~\ref{prop:betanumformula}, we
will see that SRS tiles actually perform this choice by arithmetical
means, that is, by picking points in beta-tiles that come from
$\mathbb{Z}[\beta]$. (We mention here also Barnsley's study on
\emph{fractal tops}, which provide a method to get rid of overlaps
in IFS attractors, see \cite[Chapter~4]{Barnsley:06}.)

Note that beta-tiles can be described as
\[
\mathcal{R}_\beta(x) = \mathop{\rm Lim}_{n\to\infty}
\big\{\Phi_\beta(\beta^n y) \;\big|\;
y\in\mathbb{Z}[\beta^{-1}]\cap[0,1),\,T_\beta^n(y)=x\big\}.
\]
Proposition~\ref{prop:betanumformula} shows that we have to consider
$\mathbb{Z}[\beta]$ instead of $\mathbb{Z}[\beta^{-1}]$ in this
formula to get a correspondence with SRS tiles (observe again that
in the present section $\mathbb{Z}[\beta] \cap [0,1)$ plays the same
role as $\Lambda_A$ plays in Section~\ref{cns}).

\begin{definition}[Integral beta-tile] \label{def:int}
Let $\beta$ be a Pisot number.
For $x \in \mathbb{Z}[\beta] \cap [0,1)$, the set
\[
\mathcal{S}_\beta(x) := \mathop{\rm Li}_{n\to\infty} \Phi_\beta\big(\beta^n \big(T_\beta^{-n}(x) \cap \mathbb{Z}[\beta]\big)\big),
\]
where $\mathop{\rm Li}$ is the lower Hausdorff limit defined in Section~\ref{defsrs}, is called \emph{integral $\beta$-tile}.
The tile $\mathcal{S}_\beta(0)$ is called \emph{central integral $\beta$-tile}.
\end{definition}

The difference between this definition and
Definition~\ref{def:fracbeta} is that any approximation of a tile is
given just by points $y\in\mathbb{Z}[\beta]\cap[0,1)$ with
$T_\beta^n(y)=x$, instead of considering all points in
$T_\beta^{-n}(x)$. The limitation to $\mathbb{Z}[\beta]$ is the core
of the selection process. However this implies that in general
$\mathcal{S}_\beta(x)$ is not a graph directed self-affine set. It
is obvious that
\[
\mathcal{R}_\beta(x) \supseteq \mathcal{S}_\beta(x),
\]
where equality holds if and only if $\beta$ is a unit Pisot number.

\subsection{From integral beta-tiles to SRS tiles}\label{subsec:new}
In the sequel, we will see how SRS-tiles are related to integral
beta-tiles by a linear transformation. We will show that SRS tiles
provide a decomposition of beta-tiles into disjoint pieces: the
process can be seen as selecting an integral representation in each
$p$-adic leaf of the central tile, for each prime divisor of the
norm of $\beta$. The main feature here is that this selection of an
integral representant can be performed in a dynamical way.

\begin{theorem}\label{linear}
Let $\beta$ be a Pisot number with minimal polynomial
$(x-\beta)(x^d+r_{d-1}x^{d-1}+\cdots+r_0)$ and $d=r+2s$ Galois
conjugates $\beta_1,\ldots,\beta_r\in\mathbb{R}$,
$\beta_{r+1},\ldots,\beta_{r+2s}\in\mathbb{C}\setminus\mathbb{R}$, ordered such that
$\beta_{r+1}=\overline{\beta_{r+s+1}},\,\ldots,\,\beta_{r+s}=\overline{\beta_{r+2s}}$.
Let
\[
x^d+r_{d-1}x^{d-1}+\cdots+r_0 =
(x-\beta_j)(x^{d-1}+q_{d-2}^{(j)}x^{d-2}+\cdots+q_0^{(j)})
\]
for $1\le j\le d$ and
\[
U=\left(\begin{array}{ccccc}
q^{(1)}_{0} & q^{(1)}_{1} & \cdots & q^{(1)}_{d-2} & 1 \\
\vdots & \vdots & &\vdots & \vdots \\
q^{(r)}_{0} & q^{(r)}_{1} & \cdots & q^{(r)}_{d-2} & 1 \\
\Re(q^{(r+1)}_{0}) & \Re(q^{(r+1)}_{1}) & \cdots & \Re(q^{(r+1)}_{d-2}) & 1 \\
\Im(q^{(r+1)}_{0}) & \Im(q^{(r+1)}_{1}) & \cdots & \Im(q^{(r+1)}_{d-2}) & 0 \\
\vdots & \vdots & &\vdots & \vdots \\
\Re(q^{(r+s)}_{0}) & \Re(q^{(r+s)}_{1}) & \cdots & \Re(q^{(r+s)}_{d-2}) & 1 \\
\Im(q^{(r+s)}_{0}) & \Im(q^{(r+s)}_{1}) & \cdots & \Im(q^{(r+s)}_{d-2}) & 0 \\
\end{array}\right) \in \mathbb{R}^{d \times d}.
\]
Then we have
\[
\mathcal{S}_\beta(\{\mathbf{r}\mathbf{x}\}) = U (M_\mathbf{r}-\beta I_d) \mathcal{T}_\mathbf{r}(\mathbf{x})
\]
for every $\mathbf{x} \in \mathbb{Z}^d$, where $\mathbf{r}=(r_0,\ldots,r_{d-1})$ and $I_d$ is the $d$-dimensional identity matrix.
\end{theorem}

\begin{proof}
Let $\mathbf{t} \in \mathcal{T}_\mathbf{r}(\mathbf{x})$ and $\mathbf{z}_{-n}$, $v_{-n}$, $n\in\mathbb{N}$, as in Proposition~\ref{prop:rep}, \emph{i.e.},
\[
\mathbf{t} = \mathbf{x} -
\sum_{n=0}^\infty M_\mathbf{r}^n(0,\ldots,0,v_{-n})^t \quad\mbox{and}\quad v_{-n+1}=\{\mathbf{r}\mathbf{z}_{-n}\} \ \mbox{for all }n\in\mathbb{N}.
\]
Set $b_{-n}=\beta v_{-n}-v_{-n+1}$ for $n\in\mathbb{N}$.
Then we obtain, by using $\tau_\mathbf{r}(\mathbf{x}) = M_\mathbf{r}\mathbf{x} + (0,\ldots,0,\{\mathbf{r}\mathbf{x}\})$ and $v_1=\{\mathbf{r}\mathbf{x}\}$,
\begin{equation} \label{eq:M-beta}
(M_\mathbf{r}-\beta I_d)\mathbf{t} = \tau_\mathbf{r}(\mathbf{x}) - \beta\mathbf{x} + \sum_{n=0}^\infty M_\mathbf{r}^n(0,\ldots,0,b_{-n})^t.
\end{equation}
Similarly to (\ref{eq:rMa}), we see that $(q_0^{(j)},\ldots,q_{d-2}^{(j)},1)$ is a left eigenvector of $M_\mathbf{r}$, in particular,
\[
(q_0^{(j)},\ldots,q_{d-2}^{(j)},1)M_\mathbf{r} = \beta_j(q_0^{(j)},\ldots,q_{d-2}^{(j)},1) \quad \mbox{for }1\le j\le d.
\]
By using (\ref{eq:M-beta}), we obtain
\begin{equation} \label{eq:qM}
(q_0^{(j)},\ldots,q_{d-2}^{(j)},1)(M_\mathbf{r}-\beta I_d)\mathbf{t} = (q_0^{(j)},\ldots,q_{d-2}^{(j)},1) \Big(\tau_\mathbf{r}(\mathbf{x}) - \beta\mathbf{x} + \sum_{n=0}^\infty\beta_j^n(0,\ldots,0,b_{-n})^t\Big).
\end{equation}
Since the minimal polynomial of $\beta$ can be decomposed as
\[
(x-\beta_j)(x^d+r_{d-1}^{(j)}x^{d-1}+\cdots+r_0^{(j)}) =
(x-\beta)(x-\beta_j)(x^{d-1}+q_{d-2}^{(j)}x^{d-2}+\cdots+q_0^{(j)}),
\]
we have
\[
r_0^{(j)}=-\beta q_0^{(j)}, \quad r_k^{(j)}=q_{k-1}^{(j)}-\beta q_k^{(j)} \mbox{ for } 1\le k\le d-2, \quad r_{d-1}^{(j)}=q_{d-2}^{(j)}-\beta,
\]
and obtain
\begin{equation} \label{eq:center}
(q_0^{(j)},\ldots,q_{d-2}^{(j)},1) (\tau_\mathbf{r}(\mathbf{x}) - \beta\mathbf{x}) = (r_0^{(j)},\ldots,r_{d-1}^{(j)})\mathbf{x} - \lfloor\mathbf{r}\mathbf{x}\rfloor = \{\mathbf{r}\mathbf{x}\}^{(j)}.
\end{equation}
Inserting (\ref{eq:center}) in (\ref{eq:qM}) yields
\begin{equation} \label{eq:tz}
\begin{array}{rl}
\displaystyle (q_0^{(j)},\ldots,q_{d-2}^{(j)},1)(M_\mathbf{r}-\beta I_d)\mathbf{t} & \displaystyle = \{\mathbf{r}\mathbf{x}\}^{(j)} + \sum_{n=0}^\infty b_{-n}\beta_j^n \\
& \displaystyle = \lim_{n\to\infty}\Big(\{\mathbf{r}\mathbf{x}\} +
\sum_{k=0}^{n-1} b_{-k}\beta^k\Big)^{(j)} =
\lim_{n\to\infty}\big(\beta^n\{\mathbf{r}\mathbf{z}_{-n}\}\big)^{(j)};
\end{array}
\end{equation}
here we used that $b_{-k}=\beta v_{-k}-v_{-k+1}$,
$v_{-n+1}=\{\mathbf{r}\mathbf{z}_{-n}\}$ and
$v_1=\{\mathbf{r}\mathbf{x}\}$. By
Proposition~\ref{prop:betanumformula}, we have
$\{\mathbf{r}\mathbf{z}_{-n}\} \in
T_\beta^{-n}(\{\mathbf{r}\mathbf{x}\})$ and
$\{\mathbf{r}\mathbf{z}_{-n}\} \in \mathbb{Z}[\beta]$, thus
\[
U(M_\mathbf{r}-\beta I_d)\mathbf{t} =
\lim_{n\to\infty}\Phi_\beta\big(\beta^n\{\mathbf{r}\mathbf{z}_{-n}\}\big) \in \mathcal{S}_\beta(\{\mathbf{r}\mathbf{x}\}).
\]

Now, let $\mathbf{u}\in \mathcal{S}_\beta(\{\mathbf{r}\mathbf{x}\})$.
Then there exists a sequence $(\mathbf{z}_{-n})_{n\in\mathbb{N}}$  such that $\{\mathbf{r}\mathbf{z}_{-n}\} \in T_\beta^{-n}(\{\mathbf{r}\mathbf{x}\})\cap \mathbb{Z}[\beta]$ and $\lim_{n\to\infty} \Phi_\beta(\beta^n \{\mathbf{r}\mathbf{z}_{-n}\}) = \mathbf{u}$.
Similarly to Proposition~\ref{prop:rep}, we can choose the sequence $(\mathbf{z}_{-n})_{n\in\mathbb{N}}$ such that $T_\beta(\{\mathbf{r}\mathbf{z}_{-n}\}) = \{\mathbf{r}\mathbf{z}_{-n+1}\}$ for all $n\ge 1$.
Set $b_{-n}=\beta\{\mathbf{r}\mathbf{z}_{-n-1}\}-\{\mathbf{r}\mathbf{z}_{-n}\}$.
Then (\ref{eq:tz}) implies that $\mathbf{u} \in U(M_\mathbf{r}-\beta I_d)\mathcal{T}_\mathbf{r}(\mathbf{x})$.
\end{proof}

\begin{remark}
We would like to emphasize that $U (M_\mathbf{r}-\beta I_d) \mathcal{T}_\mathbf{r}(\mathbf{x}) = \mathcal{S}_\beta(\{\mathbf{r}\mathbf{x}\})$ does not imply that the ``center'' $\mathbf{x}$ of $\mathcal{T}_\mathbf{r}(\mathbf{x})$ is mapped to the ``center'' $\Phi_\beta(\{\mathbf{r}\mathbf{x}\})$ of $\mathcal{S}_\beta(\{\mathbf{r}\mathbf{x}\})$.
Indeed, by (\ref{eq:center}), we have
\[
\Phi_\beta(\{\mathbf{r}\mathbf{x}\}) = U (\tau_\mathbf{r}(\mathbf{x}) - \beta\mathbf{x}) = U (M_\mathbf{r} - \beta I_d) \mathbf{x} + U (0,\ldots,0,\{\mathbf{r}\mathbf{x}\})^t.
\]
In particular, this means that, even though there is only a finite
number of shapes in the unit Pisot case, no SRS tile is obtained by
a $\mathbb{Z}^d$-translation of another SRS tile.
\end{remark}

\begin{corollary}
Integral beta-tiles can be written as
\[
\mathcal{S}_\beta(x) = \mathop{\rm Lim}_{n\to\infty} \Phi_\beta\big(\beta^n \big(T_\beta^{-n}(x) \cap \mathbb{Z}[\beta]\big)\big),
\]
where $\mathop{\rm Lim}$ denotes the Hausdorff limit.
\end{corollary}

We deduce the following reformulation of the set equation in Theorem~\ref{seteq}.

\begin{corollary}\label{betaseteq}
For a Pisot number $\beta$ and $x \in \mathbb{Z}[\beta]\cap[0,1)$, we have
\[
\mathcal{S}_\beta(x) = \bigcup_{y\in T_\beta^{-1}(x)\cap\mathbb{Z}[\beta]} \Lambda_\beta \mathcal{S}_\beta(y),
\]
where
$\Lambda_\beta=\mathrm{diag}\Big(\beta_1,\ldots,\beta_r,
{\Big(\begin{array}{cc} \Re(\beta_{r+1}) & \Im(\beta_{r+1}) \\
-\Im(\beta_{r+1}) & \Re(\beta_{r+1})\end{array}\Big)}, \ldots,
{\Big(\begin{array}{cc} \Re(\beta_{r+s}) & \Im(\beta_{r+s}) \\
-\Im(\beta_{r+s}) & \Re(\beta_{r+s})\end{array}\Big)}\Big)$.
\end{corollary}

The main result of this section extends the results in \cite{Akiyama:02,Ito-Rao:06,Kalle-Steiner:09} on tiling properties for unit Pisot numbers to arbitrary Pisot numbers.

\begin{theorem}
Let $\beta$ be a Pisot number.
Then the following assertions hold.
\begin{itemize}
\item
The collection $\{\mathcal{S}_\beta(x) \mid x\in\mathbb{Z}[\beta]\}$ forms a weak $m$-tiling of $\mathbb{R}^d$ for some $m\ge1$.
\item
If $\beta$ satisfies the finiteness property (F), then $\{\mathcal{S}_\beta(x) \mid x\in\mathbb{Z}[\beta]\}$ forms a weak tiling of~$\mathbb{R}^d$.
\end{itemize}
\end{theorem}

\begin{proof}
By Theorem~\ref{linear}, it is equivalent to consider the collection $\{\mathcal{T}_\mathbf{r}(\mathbf{z}) \mid \mathbf{z}\in\mathbb{Z}^d\}$ for $\mathbf{r}=(r_0,\ldots,r_{d-1})$ such that $(x-\beta)(x^d+r_{d-1}x^{d-1}+\cdots+r_0)$ is the minimal polynomial of $\beta$.
In view of Proposition~\ref{prop:beta}, the first assertion follows from Theorem~\ref{theo:mdense}, while the second assertion is a consequence of Corollary~\ref{cor:weaktiling}.
\end{proof}

\begin{example}
Consider the polynomial $x^3-3x^2+1$.
Its largest root in modulus is $\beta\approx 2.879$, a unit Pisot number.
Set $\mathbf{r}=(r_0,r_1)$ with $r_0=-1/\beta$, $r_1=-1/\beta^2$.
The $25$ SRS-tiles $\mathcal{T}_\mathbf{r}(\mathbf{x}) \subset \mathbb{R}^2$ with $\|\mathbf{x}\|_\infty \leq 2$ are shown in Figure~\ref{t30m1} on the left.

\begin{figure}[ht]
\centering
\begin{minipage}{.43\textwidth}
\vspace{2.905cm}
\includegraphics[width=\textwidth]{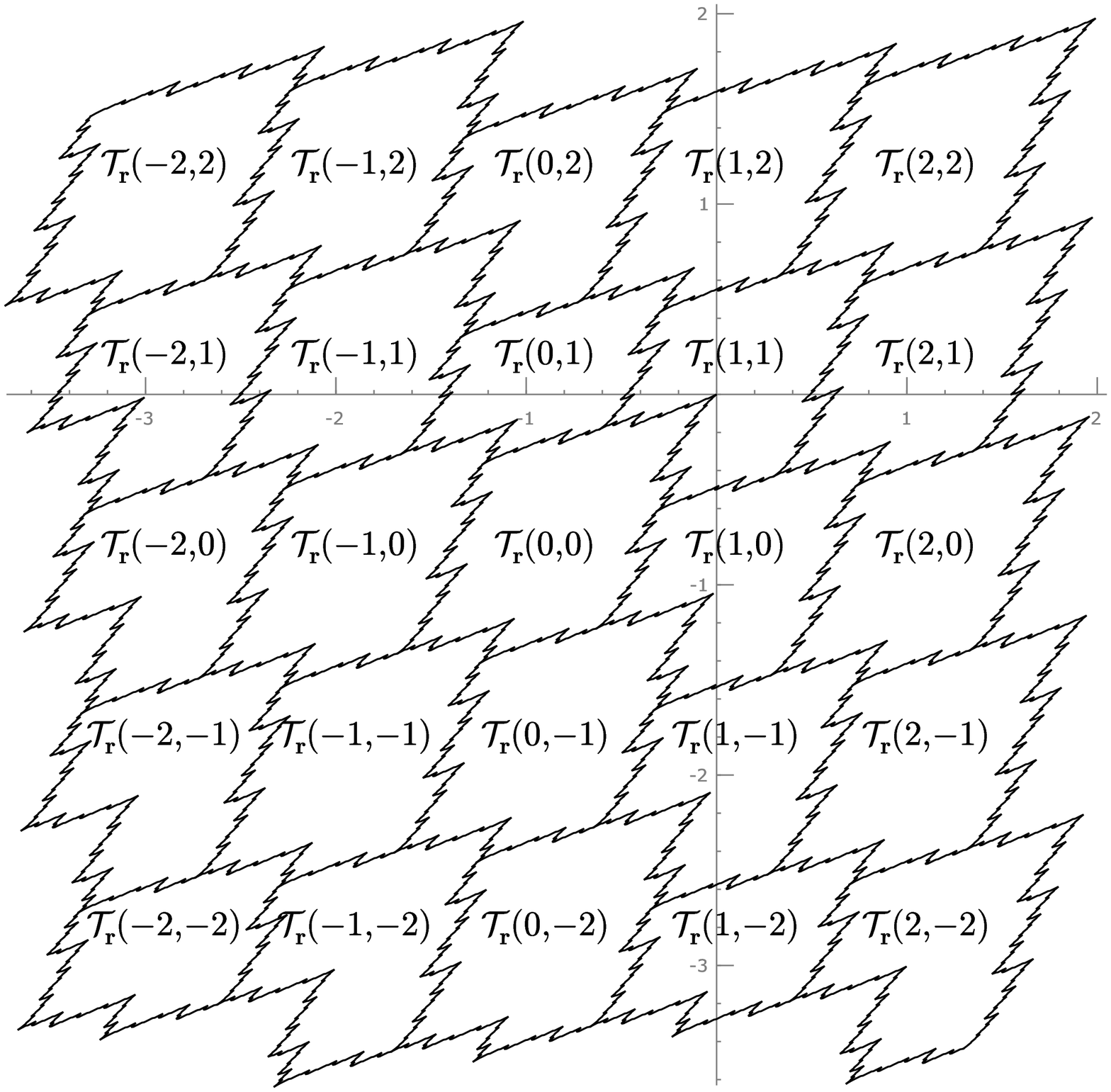}
\end{minipage}
\quad
\begin{minipage}{.53\textwidth}
\includegraphics[width=\textwidth]{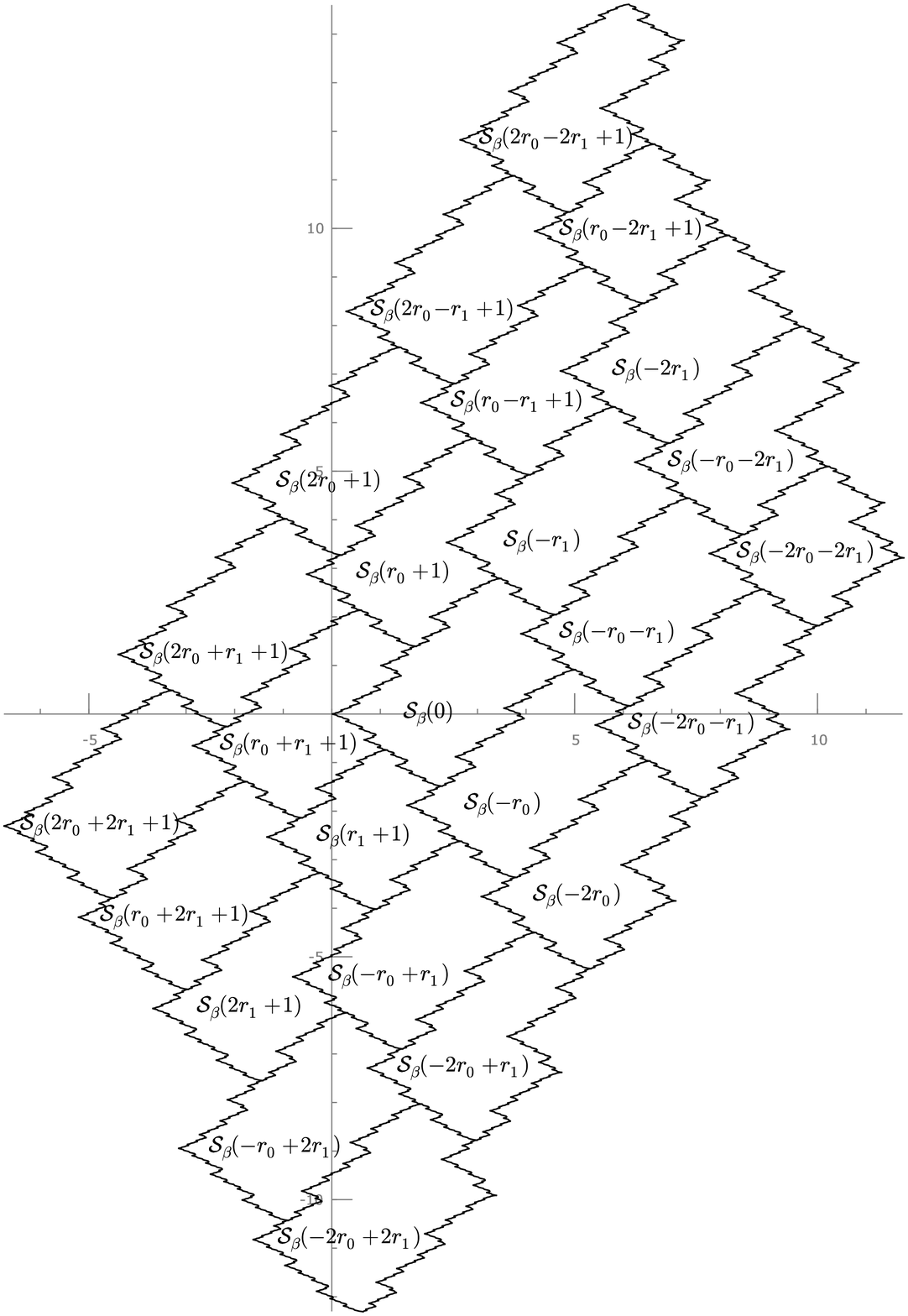}
\end{minipage}
\caption{SRS tiles $\mathcal{T}_\mathbf{r}(\mathbf{x})$ associated with $\mathbf{r}=(-1/\beta,-1/\beta^2)$, $\beta^3=3\beta^2-1$, (left) and the corresponding $\beta$-tiles $\mathcal{R}_\beta(\{\mathbf{r}\mathbf{x}\})=\mathcal{S}_\beta(\{\mathbf{r}\mathbf{x}\})$ (right).}\label{t30m1}
\end{figure}

Let $\beta_1, \beta_2$ be the Galois conjugates of $\beta$. The
numbers $\beta_1$ and $\beta_2$ are real numbers and thus the (integral) $\beta$-tiles $\mathcal{R}_\beta(\{\mathbf{r}\mathbf{x}\})=\mathcal{S}_\beta(\{\mathbf{r}\mathbf{x}\})$ are obtained by multiplying
$\mathcal{T}_\mathbf{r}(\mathbf{x})$ by the matrix
\[
U(M_\mathbf{r}-\beta I_2) = {\left(\begin{array}{cc} -\beta_2 & 1
\\ -\beta_1 & 1\end{array}\right)}
\left({\left(\begin{array}{cc} 0 & 1 \\ -r_0 & -r_1
\end{array}\right)} - {\left(\begin{array}{cc} \beta & 0 \\ 0
& \beta\end{array}\right)}\right).
\]
They are shown on the right hand side of Figure~\ref{t30m1}. Note that $\beta$ does not satisfy (F). Nevertheless, we obtain a tiling in this case.
\end{example}

\section{Conjectures and perspectives}

We are convinced that there is some potential in the study of SRS tiles that will help to gain new insights in the arithmetic as well as the geometric structure of generalized number systems. In the present section we shall discuss some possible directions of future research and state some conjectures related to SRS tiles.

In Siegel~\cite{Siegel:02}, tiles with $p$-adic factors related to non-unit Pisot substitutions have been studied. It is a natural question whether one can associate SRS type tiles to this class of substitutions. Moreover, it would be interesting to investigate the relation between Siegel's tiles and SRS type tiles. On the other hand, it should be possible to define self-similar tiles with $p$-adic factors associated with CNS with non-monic polynomials.

Another topic are geometric and topological properties of SRS tiles. As mentioned in Example~\ref{ex:pp}, there exist SRS tiles that consist of a single point. We conjecture that an SRS tile is either a single point or it is the closure of its interior. Moreover, we expect that the boundary of an SRS tile associated with a reduced parameter is not ``too big''. More precisely, we formulate the following conjecture.

\begin{conjecture}
Let $\mathbf{r} \in {\rm int}(\mathcal{D}_d)$ be a reduced parameter. Then the boundary $\partial \mathcal{T}_{\mathbf{r}}(\mathbf{x})$ has zero $d$-dimensional Lebesgue measure for each $\mathbf{x}\in\mathbb{Z}^d$.
\end{conjecture}

So far, this is only known for parameters associated with Pisot units and monic CNS because in these cases the set equation \eqref{eq:decomposition} becomes a GIFS. This leads to the following problem.

\begin{problem}\label{q1}
Characterize all parameters $\mathbf{r}\in\mathrm{int}(\mathcal{D}_d)$ for which the set equation of $\mathcal{T}_{\mathbf{r}}(\mathbf{x})$ is a (finite) GIFS equation for each $\mathbf{x}\in\mathbb{Z}^d$.
\end{problem}

Solving this problem would yield a characterization of all ``simple'' number systems related to SRS. If equation \eqref{eq:decomposition} does not give rise to a GIFS, then the structure of the SRS tiles (and the associated number systems) becomes more complicated. In \eqref{cor:515}, we proved that for $\mathbf{r}=(-2/3)$ there are infinitely many shapes of SRS tiles. A natural problem is the extension of this result to higher dimensions. For instance, is it true that the SRS tiles corresponding to a non-monic CNS have infinitely many different shapes? The ultimate goal here (which is closely related to Problem~\ref{q1}) is a characterization of all parameters $\mathbf{r}\in\mathrm{int}(\mathcal{D}_d)$ whose related SRS tiles have only finitely many different shapes.

By inspecting the pictures of SRS tiles provided throughout the paper one can see that the central tile is sometimes connected and sometimes not. Therefore we can define a Mandelbrot set
\[
\mathcal{M}_d := \{\mathbf{r}\in\mathrm{int}(\mathcal{D}_d) \mid \mathcal{T}_\mathbf{r}(\mathbf{0}) \mbox{ is connected}\}.
\]
So far, nothing is known about this set. It would be nice to find a  fast algorithm to decide whether a central tile is connected or not in order to produce an approximative picture of $\mathcal{M}_d$. After that, properties of this set can be studied.

A central feature of SRS tiles are their tiling properties. Although we could establish several results in this direction, many things remain to be done. Firstly, we were not able to prove that SRS tiles always give rise to a multiple tiling. More precisely, the following conjecture remains open.

\begin{conjecture}
Let $\mathbf{r}=(r_0,\ldots,r_{d-1}) \in
\mathrm{int}(\mathcal{D}_d)$ with $r_0 \ne 0$.
Then there is a positive integer $m$ such that $\{\mathcal{T}_{\mathbf{r}}(\mathbf{x}) \mid \mathbf{x}\in\mathbb{Z}^d\}$ is a weak $m$-tiling.
\end{conjecture}

Indeed, we conjecture that $m$ can always chosen to be equal to one, which, more precisely, reads as follows.

\begin{conjecture}\label{c74}
Let $\mathbf{r}=(r_0,\ldots,r_{d-1}) \in
\mathrm{int}(\mathcal{D}_d)$ with $r_0 \ne 0$.
Then $\{\mathcal{T}_{\mathbf{r}}(\mathbf{x}) \mid \mathbf{x}\in\mathbb{Z}^d\}$ is a weak tiling.
\end{conjecture}

Note that this implies that for all parameters associated with Pisot units the {\it Pisot conjecture} is true for beta-numeration which means that each Pisot number $\beta$ of degree $d$ gives rise to a tiling of the $(d-1)$-dimensional real vector space (see \cite{Akiyama:02}). Proving Conjecture~\ref{c74} would be a big step towards the proof of the {\it Pisot conjecture for unit Pisot substitutions} (see \cite{Arnoux-Ito:01,Ito-Rao:06}).

\bibliographystyle{siam}
\bibliography{paulibib}
\end{document}